\documentclass[12pt]{amsart}
\usepackage{amssymb,amsmath,amsthm}
\usepackage{epsfig}

\newtheorem{theorem}{Theorem}[section]
\newtheorem{lemma}[theorem]{Lemma}

\newtheorem{corollary}[theorem]{Corollary}

\newtheorem{definition}[theorem]{Definition}
\newtheorem{example}[theorem]{Example}
\newtheorem{construction}[theorem]{~~}
\newtheorem{conjecture}[theorem]{Conjecture}

\begin{document}
\pagestyle{myheadings} \thispagestyle{empty} \markboth {{\sc
 J.H.Przytycki}} {{\sc $t_k$ moves links}}
\begin{center}
\begin{LARGE}
\baselineskip=10pt {\bf  $t_k$ moves on links}

\end{LARGE}

 \ J\'OZEF H.~PRZYTYCKI
\end{center}
\begin{quotation}
{\bf Abstract: }It is a natural question to ask whether two links are 
equivalent by the following moves 
{\psfig{figure=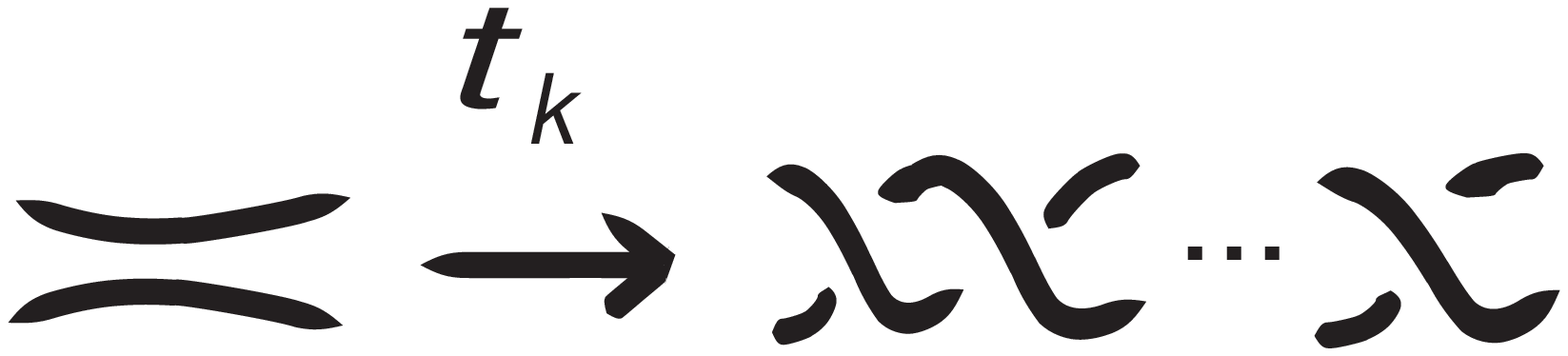,height=0.8cm}} 
(where $k$ is a fixed number of positive half twists) and if they are, 
how many moves are needed to go from one link to the other. 
In particular if $k=2$ and the second link is a trivial link 
it is the question about the unknotting number. 
The new polynomial invariants of links  often allow us to answer 
the above questions. Also the first homology groups of 
cyclic branch covers over links provide some interesting information. 
\end{quotation}

\vspace{30pt} 

 {\bf Introduction.} In the first part of the paper we apply the Jones-Conway (Homfly) and Kauffman polynomials to find whether two links are not $t_k$ equivalent and if they are, to gain some information how many moves are needed to go from one link to the other.
\par
In the second part we describe the Fox congruence classes and their relations with $t_k$ moves. We use the Fox method to analyse relations between $t_k$ moves and the first homology groups of  branched cyclic covers of links.
\par
In the third part we consider the influence of $t_k$ moves on the Goeritz and Seifert matrices and analyse Lickorish-Millett \cite{L-M-2} and Murakami \cite{Mur-1,Mur-2} formulas from the point of view of $t_k$ moves and illustrate them by various examples. At the end of the paper we outline some relations with signatures of links and non-cyclic coverings of link spaces.
\par
Now we will formulate the basic definitions and state the main 
results of the paper concerning connections between $t_k$ moves 
and the Jones-Conway polynomial invariants of links.

Consider diagrams of oriented links $L_0$ and $L_k$ which are identical, except the parts of the diagrams shown on Fig. 0.1.

\vspace*{0.8in} \centerline{\psfig{figure=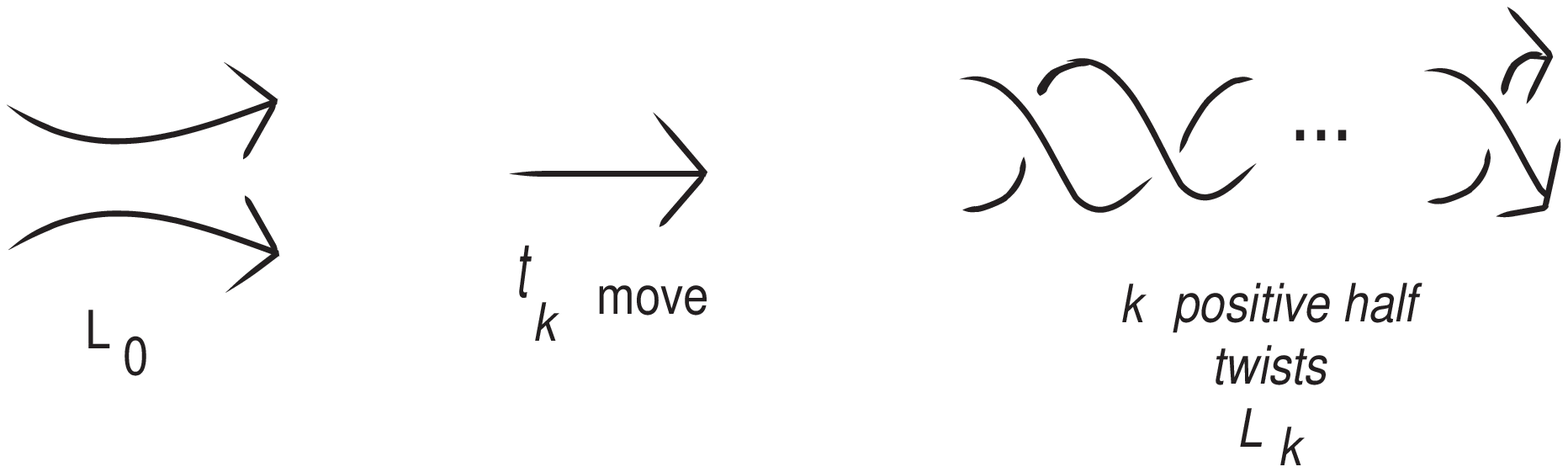,height=4.5cm}}
\begin{center}
Fig. 0.1
\end{center}

\begin{definition}
The $t_k$ move (or $k$ twist) is the elementary operation on an oriented diagram $L_0$ resulting in $L_k$ (Fig. 0.1). Two 
oriented links $L$ and $L'$ are said to be $t_k$ equivalent ($L\sim_{t_k} L'$) if one can go from $L$ to $L'$ using $t_{k}^{\mp 1}$ moves (and isotopy).
 The $t_k$ distance between $t_k$ equivalent links $L$ and $L'$ (denoted $|L,L'|_{t_k}$) is defined to be the minimal number of $t_{k}^{\mp 1}$ moves needed to go from $L$ to $L'$.
\par
The $t_k$ level distance between  $L$ and $L'$ (denoted $|L,L'|_{t_k}^{\rm lev}$) is defined to be the number of $t_k$ moves minus the number of $t_{k}^{-1}$ moves needed when we go from $L$ to $L'$ (we will show later (Corollary 1.2) that for $k>2$ it does not depend on the choice of a path joining $L$ and $L'$).
\end{definition}
\vspace{20pt}

The classical unknotting number is the $t_2$ distance from a given link to an unlink.
\vspace{20pt}

\begin{corollary}
Let $P_L(a,z)$ be a Jones-Conway
 polynomial described by the properties

\begin{enumerate}
\item[(i)]
$P_{T_1}(a,z)=1$,
\item[(ii)]
$aP_{\psfig{figure=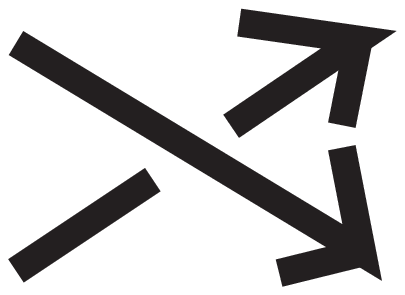,height=0.4cm
}}(a,z)+a^{-1}P_{\psfig{figure=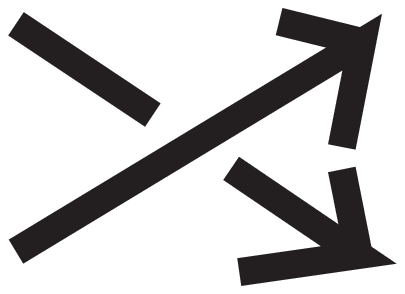,height=0.4cm
}}(a,z)=zP_{\psfig{figure=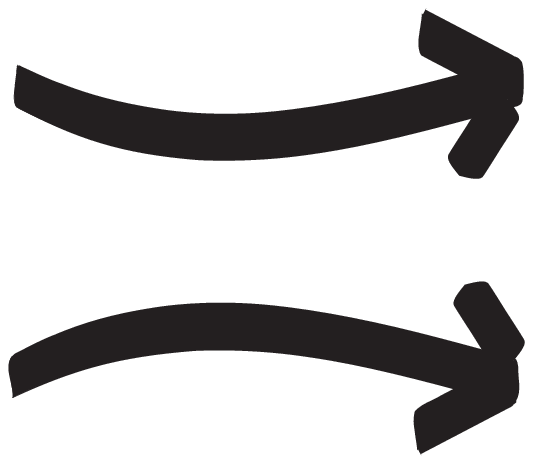,height=0.4cm }}(a,z),$
\end{enumerate}
where $T_1$ is a trivial knot. Then for 
$z_0=2\cos (\pi m/k)$ ($z_0\neq 0, \mp 2$)

$$P_{t_k(L)}(a,z_0)=(-1)^{m}a^{-k}P_L(a,z_0)$$

and for $t_k$ equivalent links $L$ and $L'$

$$P_{L'}(a,z_0)=((-1)^ma^{-k})^{|L,L'|_{t_k}^{\rm
lev}}P_L(a,z_0)$$

and neither side is identically zero.
\end{corollary}

\vspace{20pt}

We can introduce a $\bar t_k$ move and $\bar t_k$ equivalence of oriented links  ($\sim_{{\bar t}_k}$) similarly to the $t_k$ move and ($\sim_{t_k}$) 
(see Fig. 0.2).

\vspace*{0.8in} \centerline{\psfig{figure=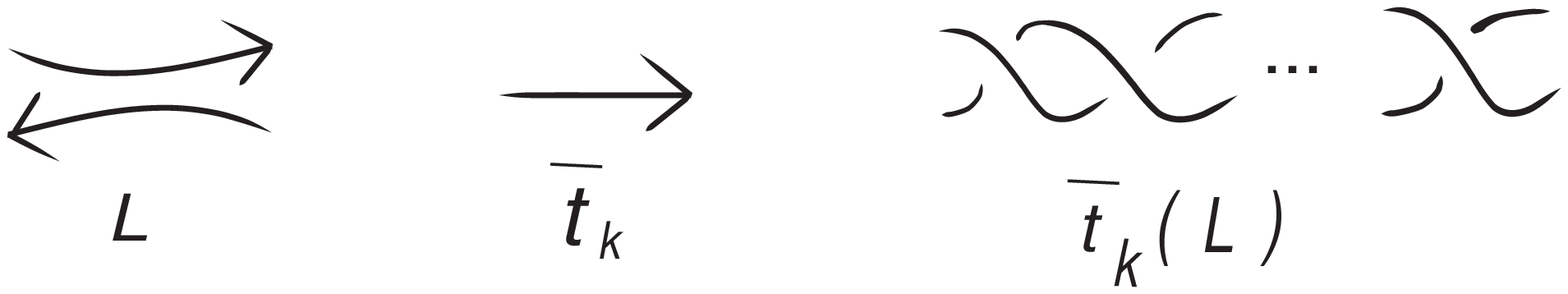,height=3.5cm}}
\begin{center}
($\bar t_k(L)$ is naturally oriented if $k$ is even)

Fig. 0.2
\end{center}

\begin{corollary}(1.8)~~
If $a_{0}^{2k}=(-1)^k$, $a_0\neq \mp i$, then 
$$P_{\bar t_{2k}(L)}(a_0,z)=P_L(a_0,z).$$
\end{corollary}

\vspace{20pt}
\begin{corollary}
Let $V_L(t)$ be the Jones polynomial described by the properties

\begin{enumerate}
\item[(i)]
$V_{T_1}(t)=1,$
\item[(ii)]
$t^{-1}V_{\psfig{figure=plus.eps,height=0.4cm}}(t)-tV_{\psfig{figure=minus.eps,height=0.4cm}}(t)=
(\sqrt{t}-\frac{1}{\sqrt{t}})V_{\psfig{figure=zero.eps,height=0.4cm}}(t),$
\end{enumerate}
then
\begin{enumerate}
\item[(a)]
If $t^k=(-1)^k$ (i.e. $t^{1/2}=-ie^{\pi i m/k}$), $t\neq -1$,
 then a $t_k$ move changes $V_L(t)$ by $(-1)^mi^k$, that is 

$$V_{t_k(L)}(t)=(-1)^mi^kV_L(t).$$

\item[(b)]
If $t^{2k}=1$ (i.e. $t=e^{\pi im/k}$), $t\neq -1$, then 
$$V_{\bar t_{2k}(L)}(t)=V_L(t).$$

\item[(c)]
Assume $k$ is odd and $t^k=-1$. Then 

$$V_{\bar t_k(L)}(t)=\omega_{4k}V_L(t),$$

where $\omega_{4k}$ is a properly chosen $4k$-root of unity (depending also on the choice of the orientation of $\bar t_k(L)$; see Theorem 1.13).
\end{enumerate}
\end{corollary}

\vspace{30pt}
\section{$t_k$-moves and Conway formulas for the Jones-Conway and Kauffman polynomials.}

When one considers the sequence of links $L, t_1(L), t_2(L), \ldots, (\psfig{figure=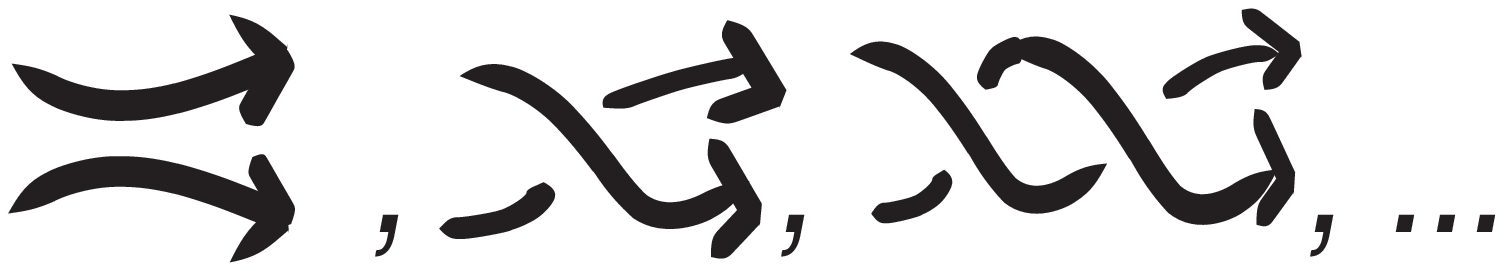,height=0.5cm} )$ then the Jones-Conway (and Kauffman) polynomials $P_L(a,z), P_{t_1(L)}(a,z), P_{t_2(L)}(a,z),\ldots $ form a (generalized) Fibonacci sequence. So one can expect that there is a nice formula which expresses $P_{t_k(L)}(a,z)$ in terms of $P_{t_1(L)}(a,z)$ and $P_{L}(a,z)$ and in fact we have the following result:

\begin{theorem}
$a^kP_{t_k(L)}(a,z)=av_{1}^{(k)}(z)P_{t_1(L)}(a,z)-v_{1}^{(k-1)}(z)P_{L}(a,z)$,
where $v_{1}^{(k+2)}(z)=zv_{1}^{(k+1)}(z)-v_{1}^{(k)}(z)$ and  
$v_{1}^{(-1)}(z)=-1$, $v_{1}^{(0)}(z)=0,$ $v_{1}^{(1)}(z)=1.$ 
In particular if one substitutes $z=p+p^{-1}$ one 
gets $v_{1}^{(k)}()=\frac{p^k-p^{-k}}{p-p^{-1}}.$ [Added for e-print: \ 
$v_{1}^{(k)}(z)$ is a variant of the Chebyshev polynomial of the 
second kind.]
\end{theorem}

\begin{proof}
We proceed by induction on $k$. For $k=1,2$ the formula from Theorem 1.1 holds:

$$aP_{t_1(L)}(a,z)=aP_{t_1(L)}(a,z)-0\cdot P_L(a,z)$$

and

$$a^2P_{t_2(L)}(a,z)=azP_{t_1(L)}(a,z)-P_L(a,z).$$

Assume that Theorem 1.1 holds for $1,2,\ldots, k-1$, ($k>2$). Now one gets:

$$a^kP_{t_k(L)}(a,z)=a^{k-1}zP_{t_{k-1}(L)}(a,z)-a^{k-2}zP_{t_{k-2}(L)}(a,z)=$$

$$=z(av_{1}^{(k-1)}(z)P_{t_1(L)}(a,z)-v_{1}^{(k-2)}(z)P_L(a,z))-$$

$$(av_{1}^{(k-2)}(z)P_{t_1(L)}(a,z)-v_{1}^{(k-3)}(z)P_L(a,z))=$$

$$a(zv_{1}^{(k-1)}(z)P_{t_1(L)}(a,z)-v_{1}^{(k-2)}(z)P_{t_1(L)}(a,z))-$$

$$(zv_{1}^{(k-2)}(z)P_L(a,z)-v_{1}^{(k-3)}(z)P_L(a,z))=$$

$$av_{1}^{(k)}(z)P_{t_1(L)}(a,z)-v_{1}^{(k-1)}(z)P_L(a,z)).$$

To see, that for $z=p+p^{-1}$, 
$v_{1}^{(k)}(z)=\frac{p^k-p^{-k}}{p-p^{-1}}$ it is enough to observe that 

$$\frac{p^{k+2}-p^{-(k+2)}}{p-p^{-1}}=(p+p^{-1})\frac{p^{k+1}-p^{-(k+1)}}{p-p^{-1}}
-\frac{p^k-p^{-k}}{p-p^{-1}}.$$

\end{proof}

\vspace{20pt}

\begin{corollary}
If $p_{0}^{2k}=1$ (i.e. $p_0=e^{\pi im/k}$), $p_0\neq \mp 1,
\mp i$ or equivalently $z_0=2\cos (\pi m/k)$; $z_0\neq 0, \mp
2$ then

$$P_{t_k(L)}(a,z_0)=(-1)^ma^{-k}P_L(a,z_0)$$

and for $t_k$ equivalent links $L$ and $L'$

$P_{L'}(a,z_0)=((-1)^ma^{-k})^{|L,L'|^{\rm lev}_{t_k}}P_L(a,z_0)$
and neither side is identically zero.
\end{corollary}

\begin{proof}
Assume $p_0\neq\mp 1, \mp i$. Then $v_{1}^{(k)}(z_0)=0$ reduces to $p_{0}^{2k}=1$, so $p_0=e^{\pi im/k}$ and $z_0=p_{0}+p_{0}^{-1}=2\cos (\pi m/k)$.
Now the equation from Theorem 1.1 reduces to 

$$a^kP_{t_k(L)}(a,z_0)=-v_{1}^{(k-1)}(z_0)P_{L}(a,z_0)=-\frac{p_{0}^{k-1}-p_{0}^{1-k}}
{p_0-p_{0}^{-1}}P_L(a,z_0)=$$

$$=p_{0}^{k}P_L(a,z_0)=(-1)^mP_L(a,z_0).$$ 
So the first part of Corollary 1.2 is proven.

For the second part it is enough to show that for each link $L$ and any 
complex number $z_0$ ($z_0 \ne 0$) $P_L(a,z_0)$ is never identically zero.
It follows from the fact that $P_L(a,a+a^{-1})\equiv 1$ (see
\cite{L-M-1}~ or~  \cite{P-1}, or apply the standard induction: it holds 
for trivial links and whenever it holds for $L_{-}(\psfig{figure=minus.eps,height=0.4cm} )$ and 
$L_{0}(\psfig{figure=zero.eps,height=0.4cm} )$ it holds for $L_{+}(\psfig{figure=plus.eps,height=0.4cm} )$ and if it holds for $L_{+}$ and $L_{0}$ it holds for $L_{-}$). 
\end{proof}

If $a=i,$ $t^{1/2}=-ip$ in the Jones-Conway polynomial $P_L(a,z)$, ($z=p+p^{-1}$), we get the (normalized) Alexander polynomial $\Delta_L(t)$ which satisfies: 
\begin{enumerate}
\item[(i)]
$\Delta_{T_1}(t)=1,$

\item[(ii)]
$\Delta_{L_{+}}(t)-\Delta_{L_{-}}(t)=(\sqrt{t}-\frac{1}{\sqrt{t}})\Delta_{L_{0}}(t).$
\end{enumerate}

\vspace{20pt}

\begin{corollary}\cite{Fo-1,Ki}
If $t^k=(-1)^k$ (i.e. $t^{1/2}=-ie^{\pi im/k}$), $t\neq -1$ then 
$\Delta_{t_k(L)}(t)=(-1)^m(-i)^k\Delta_L(t).$
\end{corollary}

\begin{proof}
It follows immediately from Corollary 1.2. One have only 
additionally notice that the formula from Corollary 1.2 remains i
true for $a=\mp i$, $p=\mp i$.
\end{proof}

\vspace{20pt}

When we substitute $a=it^{-1}$, $p=it^{1/2}$ in $P(a,z)$
($z=p+p^{-1}$) we get the Jones polynomial $V_L(t)$ which satisfies:

\begin{enumerate}
\item[(i)]
$V_{T_1}(t)=1,$
\item[(ii)]
$\frac{1}{t}V_{L_+}(t)-tV_{L_-}(t)=(\sqrt{t}-\frac{1}{\sqrt{t}})V_{L_{0}}(t).$
\end{enumerate}
There has been some confusion as to the conventions. We use that of \cite{Jo-2}.

\vspace{20pt}
\begin{corollary}
If $t^k=(-1)^k$ (i.e. $t^{1/2}=-ie^{\pi im/k}$), $t\neq -1$, 
then a $t_k$ move changes $V_L(t)$ by $(-1)^mi^k$ that is

$$V_{t_k(L)}(t)=(-1)^mi^kV_L(t).$$
\end{corollary}

\begin{proof}
It is true for $t=1$ (then $k$ is even). 
For $t\neq \mp 1$ it follows immediately from Corollary 1.2.
\end{proof}

\vspace{20pt}

\begin{corollary}
If $p_0=\varepsilon=\mp 1$ (so $z_0=2\varepsilon=\mp 2$) then
\begin{enumerate}

\item[(a)]
$a^kP_{t_k(L)}(a,z_0)=\varepsilon kaP_{t_1(L)}(a,z_0)-\varepsilon
(k-1)P_L(a,z_0)$ and

\item[(b)]$P_{t_k(L)}(a,z_0)\equiv \varepsilon a^{-k}P_L(a,z_0)({\rm
mod}~k/2^i)$ i.e. the equality holds if $P_L(a,z_0)$ is understood to be a Laurent polynomial in a with coefficients in the ring ${\bf Z}[1/2]/k{\bf Z}[1/2]$.
\end{enumerate}
\end{corollary}

\vspace{20pt}

\begin{corollary}(Generalized Conway formula).
The following formula holds for the Jones-Conway polynomial:

$$a^kP_{L_{+k}}(a,z)+a^{-k}P_{L_{-k}}(a,z)=w_{1}^{(k)}(z)P_{L_0}(a,z),$$

where $w_{1}^{(0)}=2$, $w_{1}^{(1)}=z$,
$w_{1}^{(k)}=zw_{1}^{(k-1)}-w_{1}^{(k-2)}$. After substituting $z=p+p^{-1}$ 
one gets $w_{1}^{(k)}=p^k+p^{-k}$.
\end{corollary}

\begin{proof}
From Theorem 1.1 one gets:

$$a^kP_{L_{+k}}=av_{1}^{(k)}P_{L_{+1}}-v_{1}^{(k-1)}P_{L_0}$$

and

$$a^{-k}P_{L_{-k}}=a^{-1}v_{1}^{(k)}P_{L_{-1}}-v_{1}^{(k-1)}P_{L_0}.$$

Adding these equations by sides one gets:

$$a^kP_{L_{+k}}+a^{-k}P_{L_{-k}}=v_{1}^{(k)}(aP_{L_{+1}}+a^{-1}P_{L_{-1}})-2v_{1}^{(k-1)}P_{L_0}=$$

$$=(zv_{1}^{(k)}-2v_{1}^{(k-1)})P_{L_0}.$$

Now substituting 
$w_{1}^{(k)}=zv_{1}^{(k)}-2v_{1}^{(k-1)}$ one gets the equation from Corollary 1.6 
(notice that $v_{1}^{(-1)}=-1$).
\end{proof}

\vspace{20pt}

Now we will get formulas for $\bar t_k$ moves analogous to those for $t_k$ moves.

\vspace{20pt}

\begin{theorem}
$P_{\bar
t_{2k}(L)}(a,z)=(-1)^{k}a^{2k}P_{L_{\psfig{figure=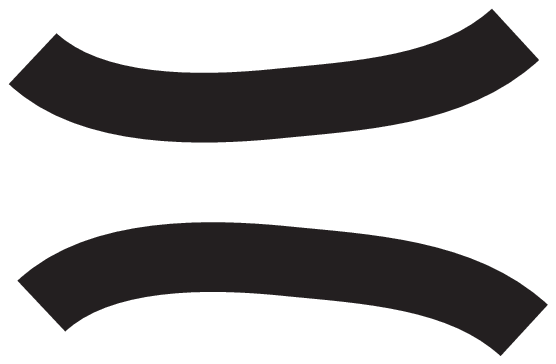,height=0.4cm}
}}(a,z)+zu_{1}^{(2k)}(a)P_{L_{\psfig{figure=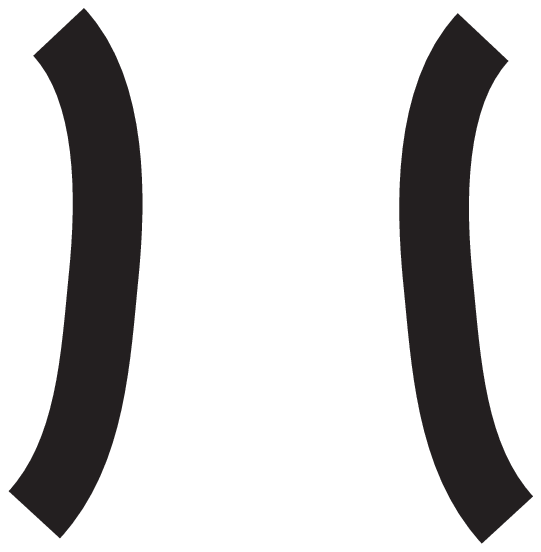,height=0.4cm}
}}(a,z)$~ where ${\bar
t_{2k}(L)},L_{\psfig{figure=zerounoriented.eps,height=0.4cm} }$ and
$L_{\psfig{figure=inftyunoriented.eps,height=0.4cm} }$ are 
oriented diagrams which are identical, except the parts of the diagrams shown on Fig. 1.1, and $u_{1}^{(0)}=0,$ 
$u_{1}^{(2)}=a,$ $u_{1}^{(2k)}=-a^2u_{1}^{(2(k-1))}+a$ or 
equivalently 
$u_{1}^{(2k)}=(-1)^{k+1}a^k\frac{a^k+(-1)^{k+1}a^{-k}}{a+a^{-1}}.$
\end{theorem}

\vspace*{0.8in} \centerline{\psfig{figure=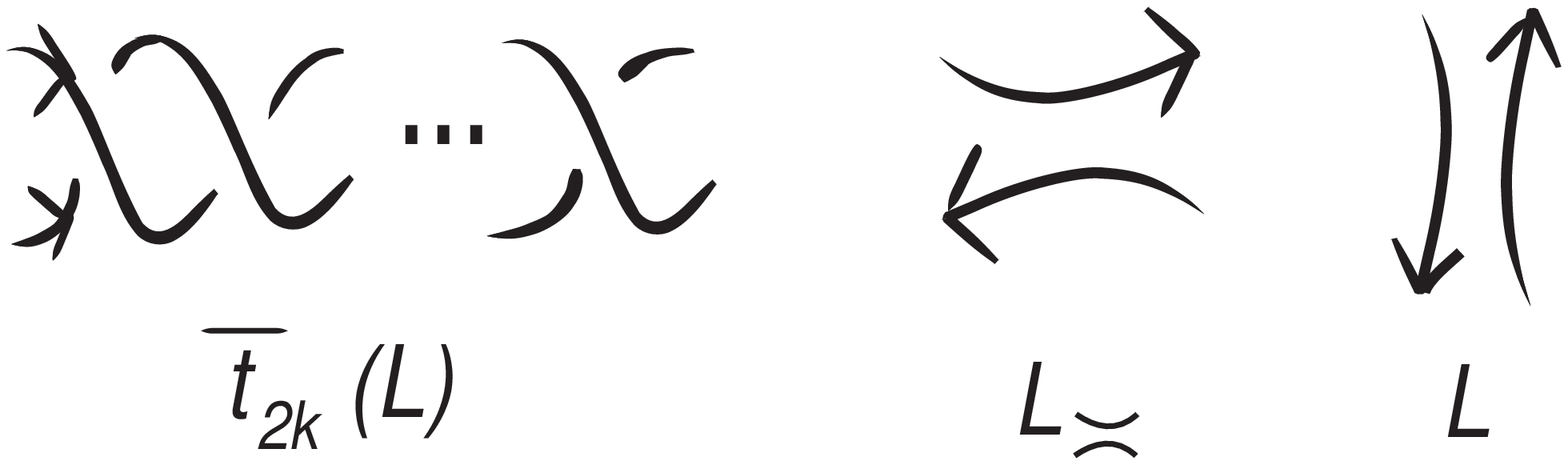,height=3cm}}
\begin{center}
Fig. 1.1
\end{center}

\begin{proof}
We proceed by induction on $k$. For $k=0,1$ the formula from Theorem 1.7 holds:

$$P_{L_{\psfig{figure=zerounoriented.eps,height=0.4cm}
}}(a,z)=P_{L_{\psfig{figure=zerounoriented.eps,height=0.4cm}
}}(a,z)+z\cdot 0\cdot
P_{L_{\psfig{figure=inftyunoriented.eps,height=0.4cm} }}(a,z)$$

and

$$P_{\bar
t_{2}(L)}(a,z)=-a^2P_{L_{\psfig{figure=zerounoriented.eps,height=0.4cm}
}}(a,z)+zaP_{L_{\psfig{figure=inftyunoriented.eps,height=0.4cm}
}}(a,z).$$

Assume that Theorem 1.7 holds for $0,1,\ldots ,k-1$ ($k\geq 2$). 
Now one gets:

$$P_{\bar t_{2k}(L)}=-a^2P_{\bar
t_{2(k-1)}(L)}+zaP_{L_{\psfig{figure=inftyunoriented.eps,height=0.4cm}
}}=-a^2((-1)^{k-1}a^{2(k-1)}P_{L_{\psfig{figure=zerounoriented.eps,height=0.4cm}}})$$

$$-a^2(zu_{1}^{(2(k-1))}P_{L_{\psfig{figure=inftyunoriented.eps,height=0.4cm}
}})+zaP_{L_{\psfig{figure=inftyunoriented.eps,height=0.4cm}
}}=(-1)^ka^{2k}P_{L_{\psfig{figure=zerounoriented.eps,height=0.4cm}
}}+$$

$$z(-a^2u_{1}^{(2(k-1))}+a)P_{L_{\psfig{figure=inftyunoriented.eps,height=0.4cm}
}}=(-1)^ka^{2k}P_{L_{\psfig{figure=zerounoriented.eps,height=0.4cm}
}}+zu_{1}^{(2k)}P_{L_{\psfig{figure=inftyunoriented.eps,height=0.4cm}
}}.$$
\end{proof}

\vspace{20pt}

\begin{corollary}
If $a_{0}^{2k}=(-1)^k$ ($a_0\neq\mp i$), then

$$P_{\bar t_{2k}(L)}(a_0,z)=P_L(a_0,z).$$
\end{corollary}

\vspace{20pt}

\begin{corollary}
If $a_0=\varepsilon i=\mp i$ then

\begin{enumerate}
\item[(a)]
$P_{\bar
t_{2k}(L)}(a_0,z)=P_{L_{\psfig{figure=zerounoriented.eps,height=0.4cm}
}}(a_0,z)+z\varepsilon
ikP_{L_{\psfig{figure=inftyunoriented.eps,height=0.4cm} }}(a_0,z)$
and
\item[(b)]
$P_{\bar t_{2k}(L)}(a_0,z)\equiv P_L(a_0,z)\ ({\rm mod\ }k)$ 
i.e. 
equality holds if $P_L(a_0,z)$ is understood to be the Laurent polynomial in $z$ with coefficients in the ring ${\bf Z}+i{\bf Z}/k({\bf Z}+i{\bf Z})$.
\item[(c)]
\cite{Fo-1} If $a=i$, $t^{1/2}=-ip$ $(z=p+p^{-1})$ one gets
 the (normalized) Alexander polynomial and  (b)
 reduces to  $\Delta_{\bar t_{2k}(L)}(t)\equiv
\Delta_L(t)\ ({\rm mod\ }k)$ i.e. the equality holds if 
$\Delta_L(t)$ is reduced to a  polynomial in ${\bf Z}_k[\sqrt{t}^{\mp 1}]$.
\end{enumerate}
\end{corollary}

\vspace{20pt}

For the Jones polynomial ($a=it^{-1}$, $p=it^{1/2}$), Corollary 1.8 reduces to:

\vspace{20pt}

\begin{corollary}
If $t^{2k}=1$, $t\neq -1$ then 

$$V_{\bar t_{2k}(L)}(t)=V_L(t).$$
\end{corollary}

\begin{proof}
It is true for $t=1$. For $t\neq \mp 1$ it follows from Corollary 1.8.
\end{proof}

\vspace{20pt}

\begin{corollary}(Generalized Conway formula).
The following formulas hold for the Jones-Conway polynomial:

\begin{enumerate}
\item[(i)]
$a^{-2k}P_{\bar t_{2k}(L)}(a,z)+a^{2k}P_{\bar
t_{2k}^{-1}(L)}(a,z)=
(-1)^k2P_{L_{\psfig{figure=zerounoriented.eps,height=0.4cm}
}}(a,z)+z(\frac{a^k+(-1)^{k+1}a^{-k}}{a+a^{-1}})^2P_{L_{\psfig{figure=inftyunoriented.eps,height=0.4cm}
}}(a,z),$

\item[(ii)]
$a^{-k}P_{\bar t_{2k}(L)}(a,z)+(-1)^ka^kP_{\bar
t_{2k}^{-1}(L)}(a,z)=((-1)^ka^k+a^{-k})P_{L_{\psfig{figure=zerounoriented.eps,height=0.4cm}
}}(a,z),$

\item[(iii)]
$a^{-k}P_{\bar t_{k}(L)}(a,z)+(-1)^{k+1}a^kP_{\bar
t_{k}^{-1}(L)}(a,z)=z(\frac{a^{-k}+\varepsilon(k)a^k}{a+a^{-1}})
P_{L_{\psfig{figure=inftyunoriented.eps,height=0.4cm} }}(a,z),$

where 

\vspace{1mm}
   \renewcommand{\arraystretch}{2}
   $\varepsilon(k)=\left\{
   \begin{array}{lr}

  -1 &
   {\rm ~if~}\ k+2 {\rm \ is~ a~ multiple~ of~} 4, \\
   1 &
   {\rm ~~~ otherwise}\ \\

   \end{array}
   \right. $
   \par
   \vspace{2mm}
\end{enumerate}
\end{corollary}

\begin{proof}
\begin{enumerate}
\item[(i)]
follows immediately  from Theorem 1.7; one has to add equations for $a^{-2k}P_{{\bar{t}}_{2k}(L)}$  and for  $a^{2k}P_{{{\bar{t}}_{2k}}^{-1}(L)}$.

\item[(ii)]
follows from Theorem 1.7, by adding equations for  $a^{-k}P_{\bar t_{2k}(L)}$ and for $(-1)^ka^{k}P_{{{\bar{t}}_{2k}}^{-1}(L)}$.

\item[(iii)]
($k$ even) follows from Theorem 1.7 by adding equations for $a^{-2k}P_{\bar t_{k}(L)}$ and for $(-1)a^{2k}P_{{\bar t_{k}}^{-1}(L)}$. If $k$ is odd then from Theorem 1.7 one gets 
$$a^{-(2k+1)}P_{\bar t_{2k+1}(L)}(a,z)=(-1)^ka^{-1}P_{\bar
t_{1}(L)}+za^{-(2k+1)}u_{1}^{(2k)}(a)P_{L_{\psfig{figure=inftyunoriented.eps,height=0.4cm}
}}(a,z)$$

and

$$a^{2k+1}P_{{\bar t_{2k+1}}^{-1}(L)}(a,z)=(-1)^kaP_{{\bar
t_{1}}^{-1}(L)}(a,z)+za^{2k+1}((-1)^{k+1}a^{-2k}u_{1}^{(2k)}(a))
P_{L_{\psfig{figure=inftyunoriented.eps,height=0.4cm} }}(a,z)$$

(in the last equality we use the fact that 
$u_{1}^{(2k)}(1/a)=(-1)^{k+1}a^{-2k}u_{1}^{(2k)}(a)$)

Adding the above equalities one gets:

$$a^{-2k+1}P_{\bar t_{2k+1}(L)}+a^{2k+1}P_{\bar
t_{2k+1}^{-1}(L)}=(-1)^k(a^{-1}P_{\bar t_1(L)}+aP_{\bar
t_{1}^{-1}(L)})+$$

$$+(z(-1)^{k-1}a^{-(k+1)}\frac{a^k+(-1)^{k+1}a^{-k}}{a+a^{-1}}+za^{k+1}
\frac{a^k+(-1)^{k+1}a^{-k}}{a+a^{-1}})P_{L_{\psfig{figure=inftyunoriented.eps,height=0.4cm}
}}=$$

$$zP_{L_{\psfig{figure=inftyunoriented.eps,height=0.4cm} }}\cdot
((-1)^k+\frac{(-1)^{k+1}a^{-1}+a^{-(2k+1)}+a^{2k+1}+(-1)^{k+1}a}{a+a^{-1}})=$$

$$=z(\frac{a^{2k+1}+a^{-(2k+1)}}{a+a^{-1}})P_{L_{\psfig{figure=inftyunoriented.eps,height=0.4cm}
}}.$$
\end{enumerate}
\end{proof}

\vspace{20pt}

 We worked, till now, with $\bar t_k$ moves for $k$ even, and the reason for this was that if $L$ is oriented then $\bar t_k(L)$ has no any natural orientation for $k$ odd. For the Jones polynomial, however, one has Jones reversing result (see \cite{L-M-2} or \cite{P-1}) so one can still find how $V_L(t)$ is changed under a $\bar t_k$ move.

Namely let $L=\{L_1,\ldots ,L_i ,\ldots ,L_n\}$ be an oriented link of $n$ components and  $L'=\{L_1,\ldots ,-L_i ,\ldots ,L_n\}$, i.e. the orientation of $L_i$ is reversed, and let $\lambda=lk(L_i,L-L_i)$.
Then 
\par
\vspace{0.1cm}\noindent
\begin{construction}
%1.12~~~~~~~~~~~~~~~~~~~~~
~~~~~~~~~~~~~~~~~~~~~$V_{L'}(t)=t^{-3\lambda}V_L(t)$.
\end{construction}
\par
%\vspace{20pt}
%\begin{theorem}
%~~
%\end{theorem}
\begin{theorem}
Consider a $\bar t_k$ move on an oriented link $L$, and assume $k$ is odd.
We have two cases:

\begin{enumerate}
\item[(i)]
$c(\bar t_k(L))<c(L)$, where  $c(L)$ denotes the number of components.  That is two components of $L$, say $L_i$ and $L_j$, are involved in the $\bar t_k$ move (see Fig.1.2). Let $\lambda=lk(L_i,L-L_i)$. 
Then for $t^k=(-1)^k$ (i.e. $t^{1/2}=-ie^{\pi im/k}$), $i\neq -1$:

$$V_{\bar
t_k(L)}(t)=(-1)^mi^kt^{-3\lambda}V_L(t)=(-1)^{m+\lambda}i^ke^{-6\pi
im\lambda/k}V_L(t),$$

where the orientation of $\bar t_k(L)$ is chosen so that it does not agree with the orientation of $L_i$.

\item[(ii)]
$c(\bar t_k(L))=c(L)$. That is one component of $L$ is involved in the $\bar t_k$ move. Let $L_{\psfig{figure=inftyunoriented.eps,height=0.4cm} }$ denote the smoothing of $L$ (Fig. 1.3).
 $L_{\psfig{figure=inftyunoriented.eps,height=0.4cm}
}$ has more components than $L$ and let $L_i$,  $L_j$ be the new components of 
$L_{\psfig{figure=inftyunoriented.eps,height=0.4cm} }$ (Fig. 1.3).
Let $\lambda=lk(L_i,L_{\psfig{figure=inftyunoriented.eps,height=0.4cm} }-L_i)$ and assume that $\bar t_k(L)$ is oriented in such a way that its orientation agrees with that of $L_{\psfig{figure=inftyunoriented.eps,height=0.4cm} }$ with exception of $L_i$. Then for $t^k=(-1)^k$ (i.e. 
$t^{1/2}=-ie^{\pi im/k}$), $i\neq -1$:

$$V_{\bar t_k(L)}(t)=t^{-3\lambda}V_L(t)=(-1)^{\lambda}e^{-6\pi
im\lambda/k}V_L(t).$$
\end{enumerate}
\end{theorem}

\begin{proof}
\begin{enumerate}
\item[(i)]
We use the Jones reversing result and Corollary 1.4 and we get (see Fig. 1.2)

$$V_{\bar
t_k(L)}(t)=(-1)^mi^kV_{L'}(t)=(-1)^mi^kt^{-3\lambda}V_L(t)=(-1)^{m+\lambda}i^ke^{6\pi
im\lambda/k}V_L(t).$$

\vspace*{0.8in} \centerline{\psfig{figure=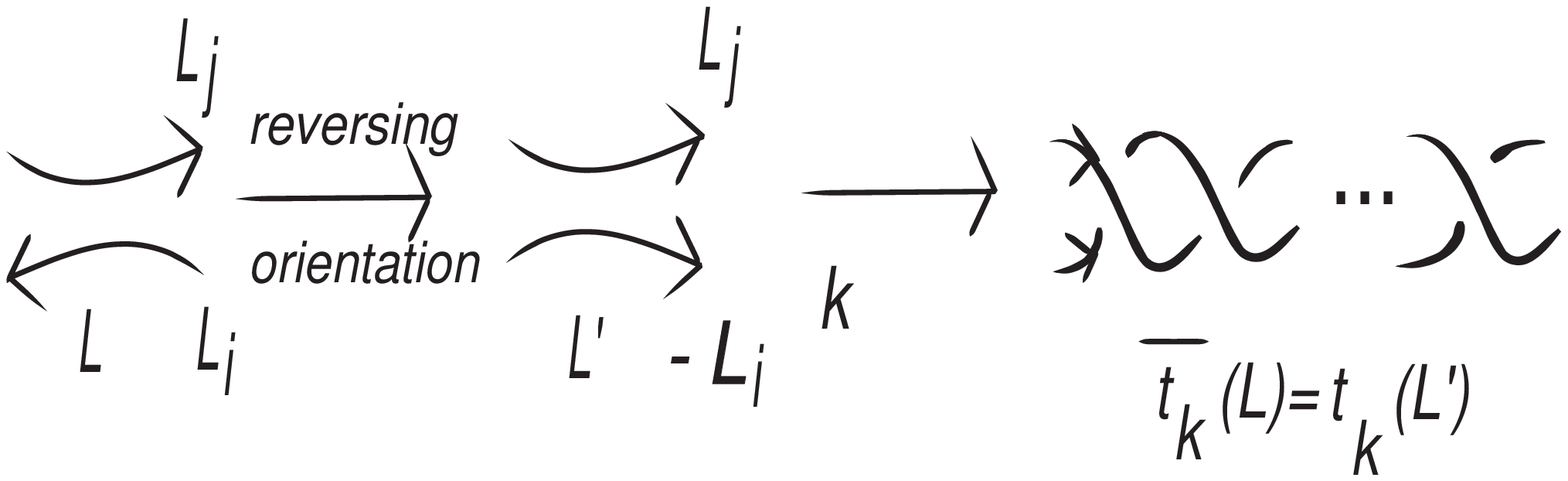,height=4.5cm}}
\begin{center}
Fig. 1.2
\end{center}

\item[(ii)]
We use Corollary 1.4 and the part (i) and we get (see Fig. 1.4):

$$V_{\bar
t_k(L)}(t)=(-1)^mi^kt^{-3\lambda+k}V_{L'}(t)=(-1)^mi^kt^{-3\lambda+k}(-1)^mi^kV_{L}(t)=
t^{-3\lambda}V_L(t)=$$

$$(-1)^{\lambda}e^{-6\pi im\lambda/k}V_L(t).$$
\end{enumerate}
\vspace{20pt} \centerline{\psfig{figure=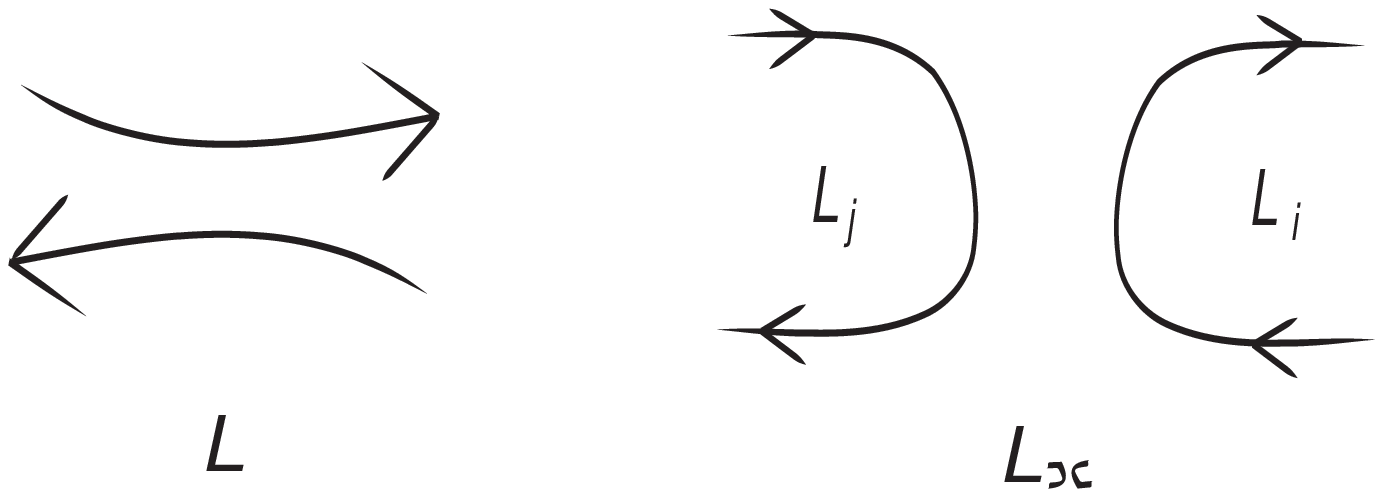,height=2.1cm}}
\begin{center}
Fig. 1.3
\end{center}

\vspace*{0.8in} \centerline{\psfig{figure=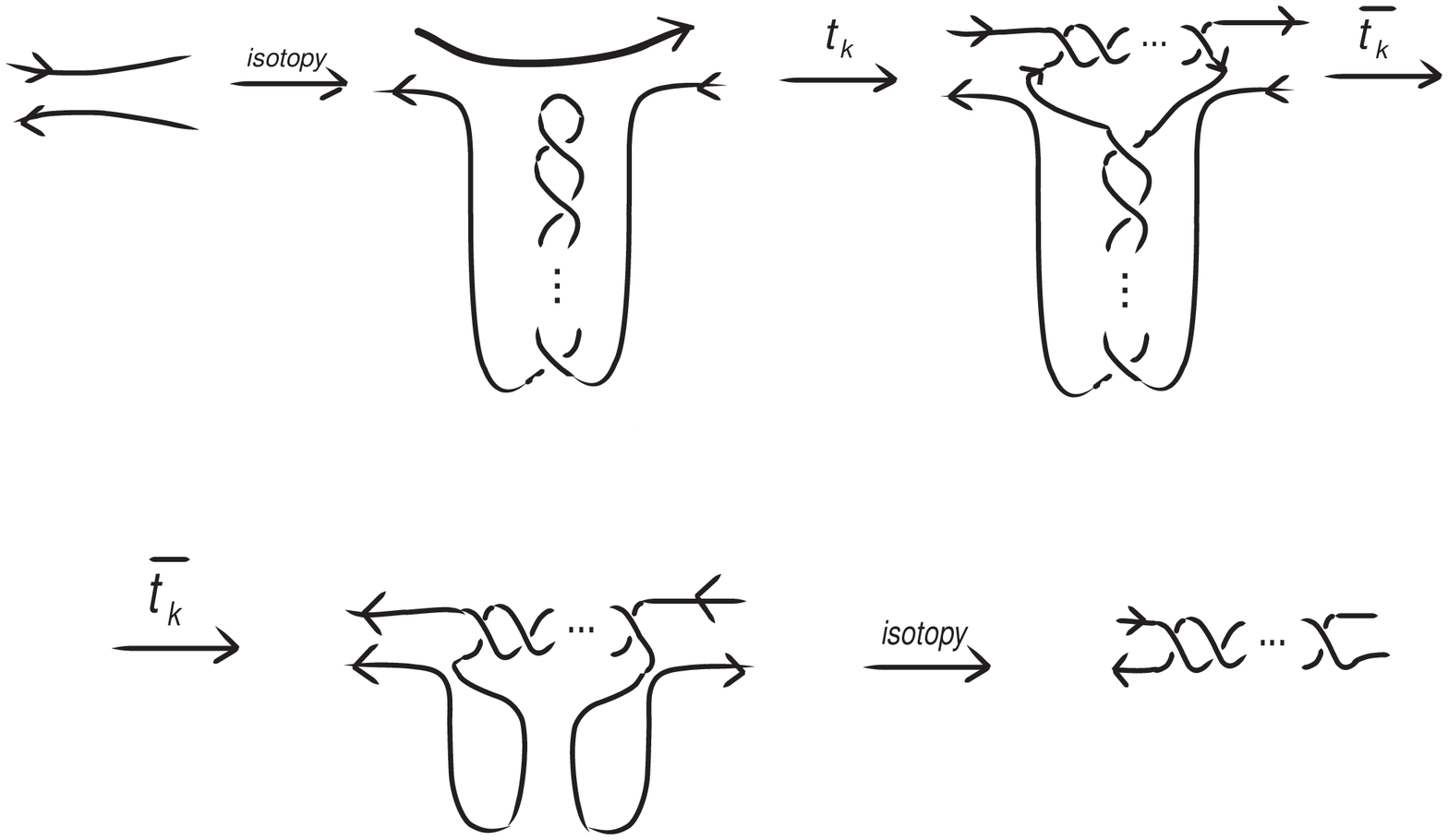,height=8.5cm}}
\begin{center}
Fig. 1.4
\end{center}
\end{proof}

\vspace{20pt}

It is possible to get Theorem 1.13 by considering the variant of the Jones polynomial which is an invariant of regular isotopy and does not depend on orientation (\cite{Ka-4}). 

We will use this idea considering how the Kauffman polynomial changes under $t_k$ and $\bar t_k$ moves.

Two  diagrams of links are regularly isotopic iff one can be obtained from 
the other by a sequence of Reidemeister moves of type $\Omega_{2}^{\mp
1}$, $\Omega_{3}^{\mp 1}$ and isotopy of the projection plane (see Fig.1.5).

\vspace*{0.8in} \centerline{\psfig{figure=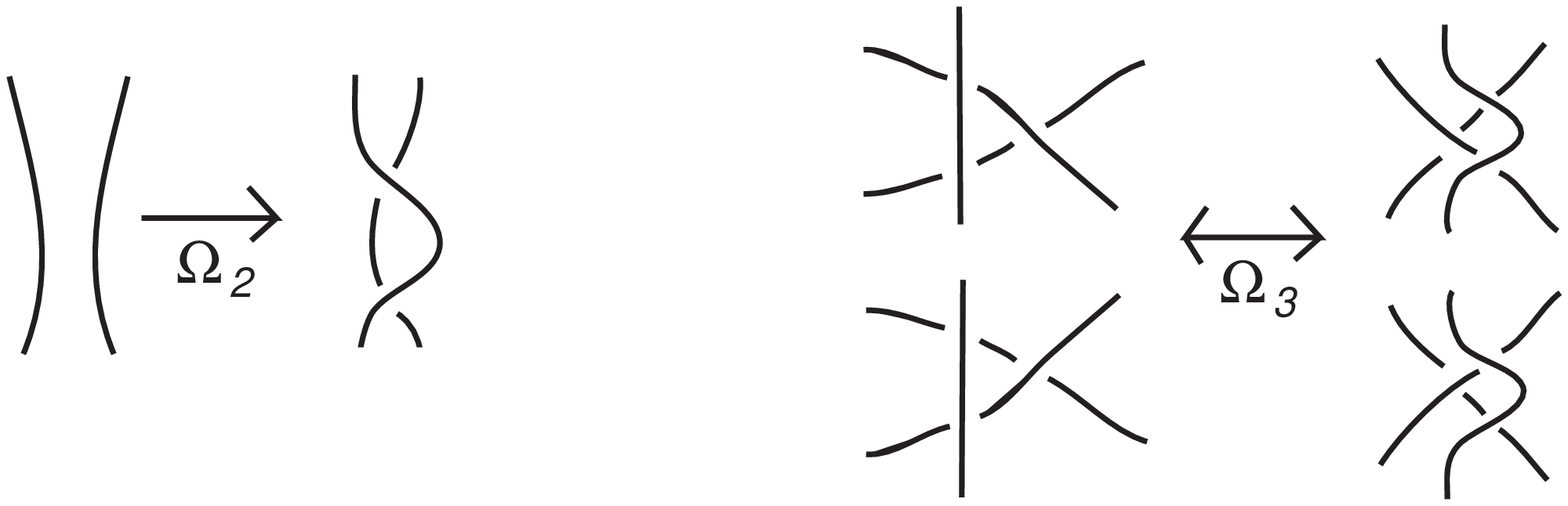,height=4.5cm}}
\begin{center}
Fig. 1.5
\end{center}

The Kauffman polynomial of regular isotopy of unoriented diagrams is defined by  (see  \cite{Ka-3}; also \cite{P-1}):

\begin{enumerate}
\item[(1)]
$\Lambda_{T_1}(a,x)=a^{{\rm tw}(T_1)}$, where $T_1$ is a diagram representing the trivial knot (up to isotopy) and ${\rm tw}(T_1)=\sum {\rm sgn\ }p$  where the sum is taken over all crossings of $T_1$.

\item[(2)]
$\Lambda_{\psfig{figure=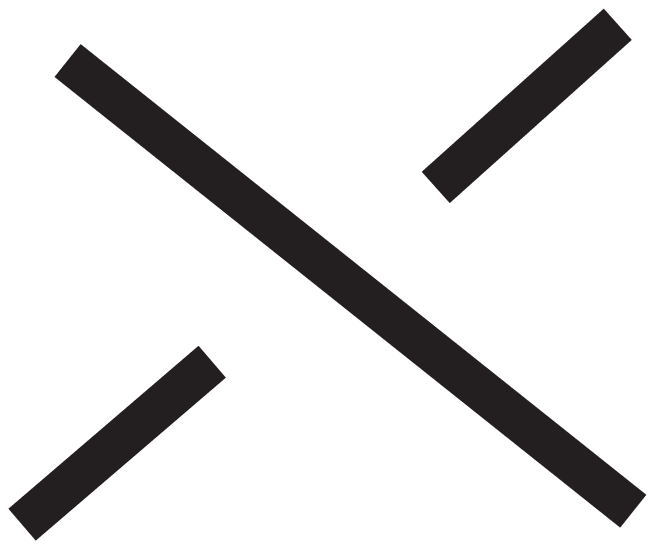,height=0.4cm}
 }(a,x)+\Lambda_{\psfig{figure=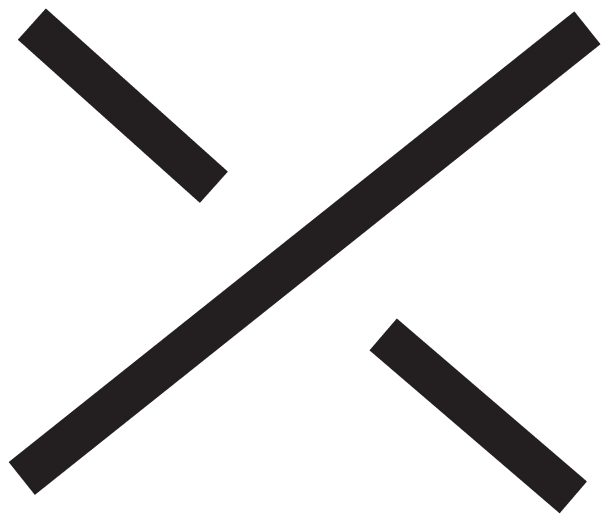,height=0.4cm}
 }(a,x)=x\Lambda_{\psfig{figure=zerounoriented.eps,height=0.4cm}
 }(a,x)+x\Lambda_{\psfig{figure=inftyunoriented.eps,height=0.4cm}
 }(a,x).$

\end{enumerate}
The Kauffman polynomial of oriented links is defined by 
$$F_L(a,x)=a^{-{\rm tw}(L)}\Lambda_{L}(a,x).$$

\vspace{20pt}

\begin{theorem}\label{Theorem 1.14}
$\underbrace{\Lambda_{\psfig{figure=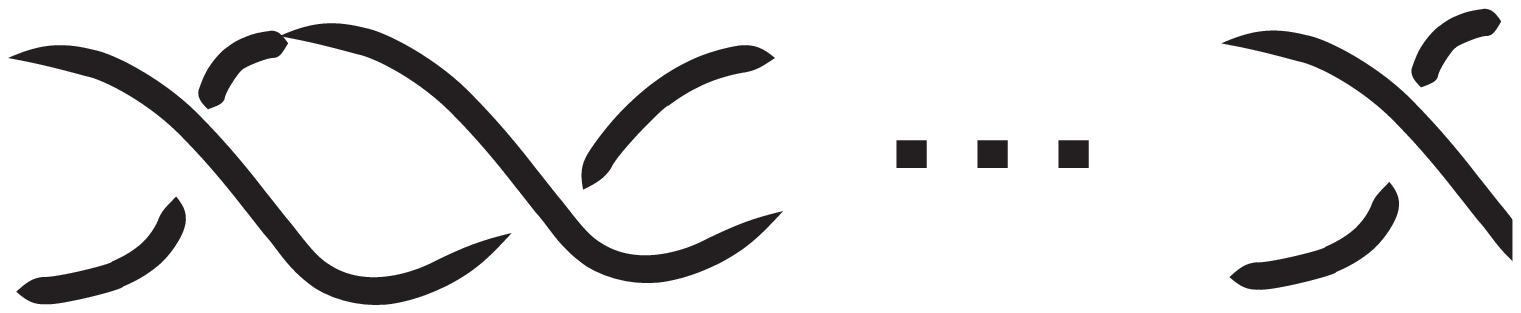,height=0.4cm}
 }}_{k {\rm\ half~twists}}(a,x)=v_{1}^{(k)}(x)\Lambda_{\psfig{figure=plusunoriented.eps,height=0.4cm}
 }(a,x)-v_{1}^{(k-1)}(x)\Lambda_{\psfig{figure=zerounoriented.eps,height=0.4cm}
 }(a,x)+xv_{2}^{(k)}(a,x)\Lambda_{\psfig{figure=inftyunoriented.eps,height=0.4cm}
 }(a,x),$
 $k$ half twists where $v_{1}^{(k)}()$ is the same as in Theorem 1.1 and  $v_{2}^{(0)}(a,x)=0$, $v_{2}^{(1)}(a,x)=0$,
 $v_{2}^{(2)}(a,x)=a^{-1}$,
 $v_{2}^{(k)}(a,x)=xv_{2}^{(k-1)}(a,x)-v_{2}^{(k-2)}(a,x)+a^{1-k}$.
 In particular for $x=p+p^{-1}$ one gets 
 $v_{1}^{(k)}()=\frac{p^k-p^{-k}}{p-p^{-1}},
 v_{2}^{(k)}=((p-p^{-1})(a+a^{-1}-(p+p^{-1}))^{-1}(-a^{-1}(p^k-p^{-k})+
 p(a^{-k}-p^{-k})-p^{-1}(a^{-k}-p^k)).$
\end{theorem}

\begin{proof}
We proceed by induction on $k$. For $k=0,1,2$ the formula from Theorem 1.14 holds. Assume that it holds for  $0,1,\ldots ,k-1$ ($k>2$). Now one gets:

$$\underbrace{\Lambda_{\psfig{figure=ktwist.eps,height=0.4cm}
 }}_{k{\rm\ half\ twists}}=\underbrace{x\Lambda_{\psfig{figure=ktwist.eps,height=0.4cm}
 }}_{(k-1){\rm\ half\ twists}}-\underbrace{\Lambda_{\psfig{figure=ktwist.eps,height=0.4cm}
 }}_{(k-2){\rm\ half\ twists}}+\ \ \ \ \ xa^{1-k}\Lambda_{\psfig{figure=inftyunoriented.eps,height=0.4cm}
 }=$$

 $$=x(v_{1}^{(k-1)}\Lambda_{\psfig{figure=plusunoriented.eps,height=0.4cm}
 }-v_{1}^{(k-1)}\Lambda_{\psfig{figure=zerounoriented.eps,height=0.4cm}
 }+xv_{2}^{(k-1)}\Lambda_{\psfig{figure=inftyunoriented.eps,height=0.4cm}
 })-(v_{1}^{(k-2)}\Lambda_{\psfig{figure=plusunoriented.eps,height=0.4cm}
 }-v_{1}^{(k-3)}\Lambda_{\psfig{figure=zerounoriented.eps,height=0.4cm}
 }+xv_{2}^{(k-2)}\Lambda_{\psfig{figure=inftyunoriented.eps,height=0.4cm}
 })+xa^{1-k}\Lambda_{\psfig{figure=inftyunoriented.eps,height=0.4cm}
 }=$$

 $$=(xv_{1}^{(k-1)}-v_{1}^{(k-2)})\Lambda_{\psfig{figure=plusunoriented.eps,height=0.4cm}
 }-(xv_{1}^{(k-2)}-v_{1}^{(k-3)})\Lambda_{\psfig{figure=zerounoriented.eps,height=0.4cm}
 }+x(xv_{2}^{(k-1)}-v_{2}^{(k-2)}+a^{1-k})\Lambda_{\psfig{figure=inftyunoriented.eps,height=0.4cm}
 }=$$

 $$=v_{1}^{(k)}\Lambda_{\psfig{figure=plusunoriented.eps,height=0.4cm}
 }-v_{1}^{(k-1)}\Lambda_{\psfig{figure=zerounoriented.eps,height=0.4cm}
 }+xv_{2}^{(k)}\Lambda_{\psfig{figure=inftyunoriented.eps,height=0.4cm}
 }.$$

 The formula for $v_{2}^{(k)}(a,p+p^{-1})$ may be verified directly but we omit this tedious task by considering the trivial links of Fig. 1.6. 
 From this figure we get immediately that

 $$a^{-k}=\frac{p^k-p^{-k}}{p-p^{-1}}a^{-1}-\frac{p^{k-1}-p^{-(k-1)}}{p-p^{-1}}
 +(p+p^{-1})v_{2}^{(k)}(a,x)\frac{a+a^{-1}-(p+p^{-1})}{p+p^{-1}},$$

 and it finishes the proof of Theorem 1.14.
\end{proof}

\vspace*{0.8in} \centerline{\psfig{figure=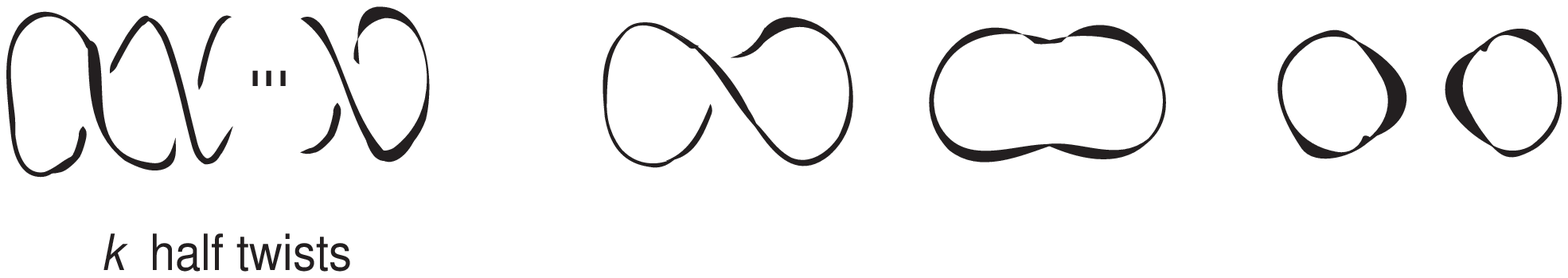,height=3cm}}
\begin{center}
Fig. 1.6
\end{center}

\vspace{20pt}

\begin{corollary}
\begin{enumerate}
\item[(a)]
If $p_{0}^{2k}=1$ (i.e. $p_{0}=e^{\pi im/k}$), $p_{0}\neq \mp
1,\ \mp i$ or equivalently $x_0=2\cos (\pi m/k)$, $x_0\neq 0,
\mp 2$, then 

$$\Lambda_{\psfig{figure=ktwist.eps,height=0.4cm}
 }(a,x_0)=(-1)^m\Lambda_{\psfig{figure=zerounoriented.eps,height=0.4cm}
 }(a,x_0)+\frac{a^{-k}-p_{0}^{-k}}{(a-p_0)(1-a^{-1}p_{0}^{-1})}\Lambda_{\psfig{figure=inftyunoriented.eps,height=0.4cm}
 }(a,x_0).$$

 \item[(b)]
 If $p_{0}^{2k}=1$, $p_{0}\neq \mp 1,\ \mp i$,
 $a_{0}^{k}=p_{0}^{k}, a_0\neq p_{0}^{\mp 1}$ then 

 $$\Lambda_{\psfig{figure=ktwist.eps,height=0.4cm}
 }(a_0,x_0)=(-1)^m\Lambda_{\psfig{figure=zerounoriented.eps,height=0.4cm}
 }(a_0,x_0).$$

 \item[(c)]
 If $p_0=\varepsilon =\mp 1$ (so $x_0=2\varepsilon =\mp 2$) and
 $a_{0}^{k}=\varepsilon^{k}$, $a_0\neq\varepsilon$ then

 $$\Lambda_{\psfig{figure=ktwist.eps,height=0.4cm}
 }(a_0,x_0)\equiv\varepsilon \Lambda_{\psfig{figure=zerounoriented.eps,height=0.4cm}
 }(a_0,x_0)({\rm mod\ }k/2^i),$$

 i.e. equality holds in the ring ${\bf Z}[a_0/2]/k{\bf Z}[a_0/2]$.
\end{enumerate}
\end{corollary}

\begin{proof}
It follows from Theorem 1.14 similarly as Corollaries 1.2 and 1.5 followed from Theorem 1.1.
\end{proof}

\vspace{20pt}

\begin{corollary}(Generalized Conway formula)
$$\underbrace{\Lambda_{\psfig{figure=ktwist.eps,height=0.4cm}
 }}_{k{\rm\ half\ twists}}(a,x)+\Lambda_{\psfig{figure=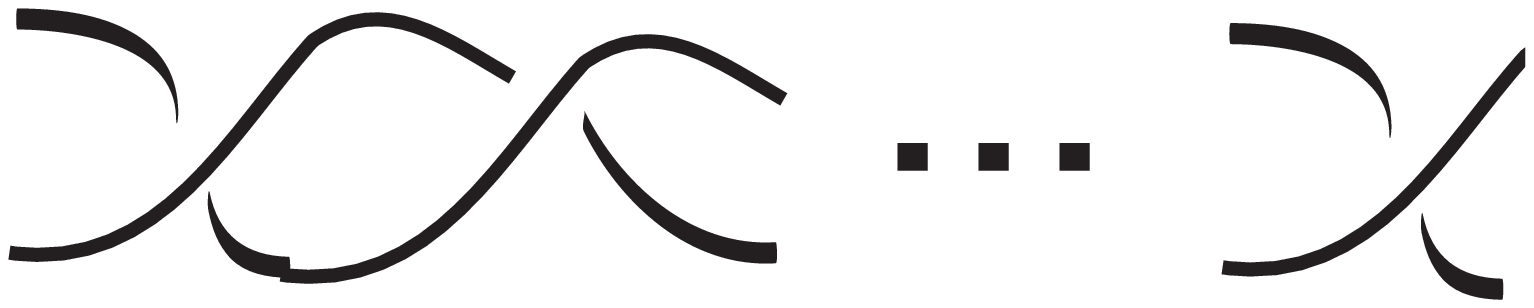,height=0.4cm}
 }(a,x)=w_{1}^{(k)}\Lambda_{\psfig{figure=zerounoriented.eps,height=0.4cm}
 }(a,x)+xw_{2}^{(k)}(a,x)\Lambda_{\psfig{figure=inftyunoriented.eps,height=0.4cm}
 }(a,x),$$

 where $w_{1}^{(k)}(x)=w_{1}^{k}(z)$ from Corollary 1.6, and 

 $$w_{2}^{(0)}(a,x)=0,\ w_{2}^{(1)}(a,x)=1,\
 w_{2}^{(k)}(a,x)=xw_{2}^{(k-1)}(a,x)-w_{2}^{(k-2)}(a,x)+a^{k-1}+a^{1-k};$$

 when one substitutes  $x=p+p^{-1}$ then 

 $w_{1}^{(k)}(x)=p^k+p^{-k}$ and 

 $$w_{2}^{(k)}(a,x)=a^{-(k-1)}p^{-(k-1)}(\frac{a^k-p^k}{a-p})(\frac{1-a^kp^k}{1-ap}).$$
\end{corollary}

\begin{proof}
From Theorem 1.14 one gets:

$$\Lambda_{\psfig{figure=ktwist.eps,height=0.4cm}
 }=v_{1}^{(k)}\Lambda_{\psfig{figure=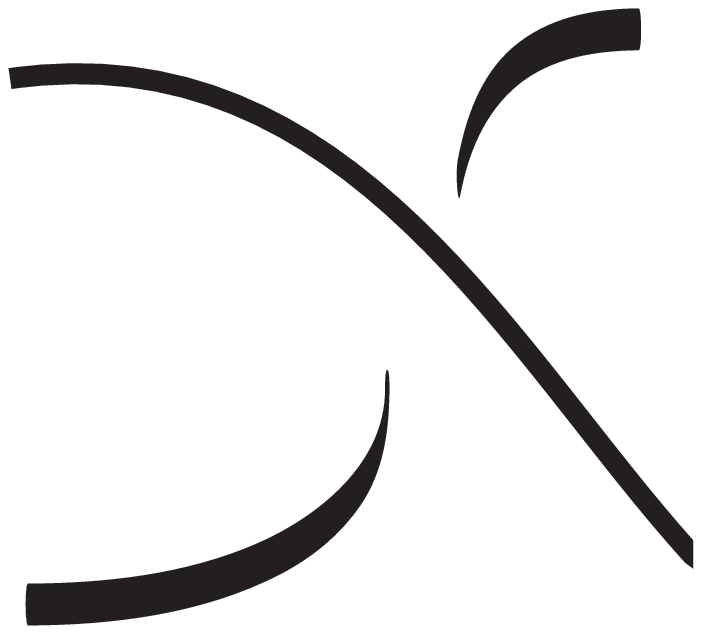,height=0.4cm}
 }-v_{1}^{(k-1)}\Lambda_{\psfig{figure=zerounoriented.eps,height=0.4cm}
 }+xv_{2}^{(k)}(a,x)\Lambda_{\psfig{figure=inftyunoriented.eps,height=0.4cm}
 }$$

and

 $$\Lambda_{\psfig{figure=minusktwist.eps,height=0.4cm}
 }=v_{1}^{(k)}\Lambda_{\psfig{figure=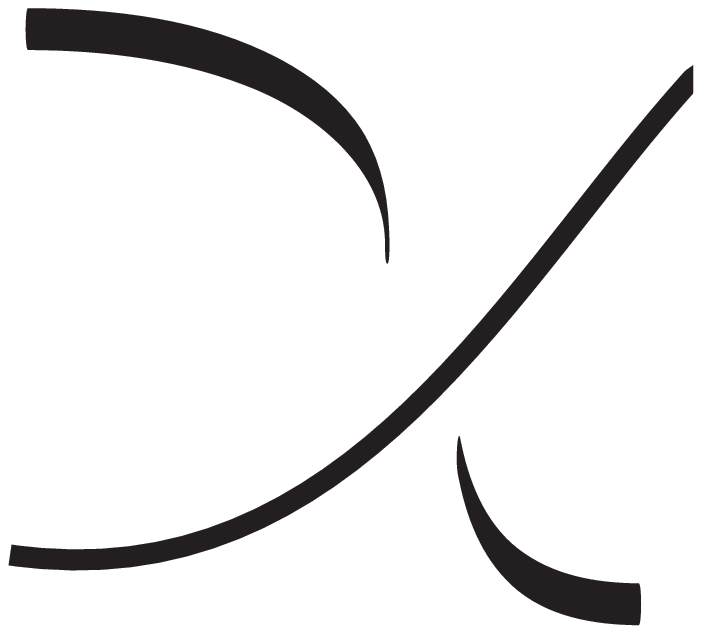,height=0.4cm}
 }-v_{1}^{(k-1)}\Lambda_{\psfig{figure=zerounoriented.eps,height=0.4cm}
 }+xv_{2}^{(k)}(a^{-1},x)\Lambda_{\psfig{figure=inftyunoriented.eps,height=0.4cm}
 }$$

 Adding the above equations by sides one gets:
 
 $$\Lambda_{\psfig{figure=ktwist.eps,height=0.4cm}
 }+\Lambda_{\psfig{figure=minusktwist.eps,height=0.4cm}
 }=v_{1}^{(k)}(\Lambda_{\psfig{figure=twist.eps,height=0.4cm}
 }+\Lambda_{\psfig{figure=minustwist.eps,height=0.4cm}
 })-2v_{1}^{(k-1)}\Lambda_{\psfig{figure=zerounoriented.eps,height=0.4cm}
 }+$$

 $$+x(v_{2}^{(k)}(a,x)+v_{2}^{(k)}(a^{-1},x))\Lambda_{\psfig{figure=inftyunoriented.eps,height=0.4cm}
 }=(xv_{1}^{k}-2v_{1}^{(k-1)})\Lambda_{\psfig{figure=zerounoriented.eps,height=0.4cm}
 }+$$

 $$+x(v_{1}^{(k)}+v_{2}^{(k)}(a,x)+v_{2}^{(k)}(a^{-1},x))\Lambda_{\psfig{figure=inftyunoriented.eps,height=0.4cm}
 }.$$

 Now substituting $w_{1}^{(k)}=xv_{1}^{k}-2v_{1}^{(k-1)}$
 and 
 $w_{2}^{(k)}=v_{1}^{k}(a)+v_{2}^{(k)}(a,x)+v_{2}^{(k)}(a^{-1},x)$
 one gets the equation from Corollary 1.16.
\end{proof}

\vspace{20pt}

We end this part of the paper by translating Corollary 1.15(b) into the Kauffman polynomial of oriented links.

\vspace{20pt}

\begin{corollary}
 If $p_{0}^{2k}=1$, $p_{0}\neq \mp 1,\ \mp i$ and 
  $a_{0}^{k}=p_{0}^{k}, a_0\neq p_{0}^{\mp 1}$ then

 \begin{enumerate}
\item[(a)]
$F_{\psfig{figure=ktwist.eps,height=0.4cm}
 }(a_0,p_0)=a_{0}^{k}a_{0}^{{\rm
 tw}(\psfig{figure=zerounoriented.eps,height=0.4cm})
 -{\rm tw}(\psfig{figure=ktwist.eps,height=0.4cm})
 }F_{\psfig{figure=zerounoriented.eps,height=0.4cm}
 }(a_0,p_0);$

 In particular

 \item[(b)]
 $F_{t_k(L)}(a_0,p_0)=F_L(a_0,p_0),$
 \item[(c)]
 $F_{\bar t_{2k}(L)}(a_0,p_0)=F_L(a_0,p_0),$
 \item[(d)]
 $F_{\bar t_{k}(L)}(a_0,p_0)=a^{4\lambda}F_L(a_0,p_0),$ where $k$ 
 is odd and $\lambda$ defined as follows (compare Theorem 1.13):

 \begin{enumerate}
\item[(i)]
If $L$ has more components than $\bar t_{k}(L)$ and 
$L_i$ is the only component of $L$ such that the chosen orientation on 
$\bar
t_{k}(L)$ does not agree with that of $L_i$ then $\lambda=lk(L_i,L-L_i)$ 
(compare Fig. 1.2).

\item[(ii)]
If $L$ has the same number of components as $\bar t_{k}(L)$, consider the smoothing 
$L_{\psfig{figure=inftyunoriented.eps,height=0.4cm}
 }$ of $L$ (Fig. 1.3). Let $L_i$ be the only component of $L$
 such that the chosen orientation on $\bar
t_{k}(L)$ does not agree with that of $L_i$ then 
$\lambda=lk(L_i,L-L_i)$
(compare Fig. 1.3).

 \end{enumerate}
 \end{enumerate}

\end{corollary}

\begin{proof}
\begin{enumerate}
\item[(a)]
follows immediately from Corollary 1.15(b) and the definition of $F_L(a,x)$.

\item[(b) and (c)] hold because in these cases 
 ${\rm tw}(\psfig{figure=zerounoriented.eps,height=0.4cm}
 )-{\rm tw}(\psfig{figure=ktwist.eps,height=0.4cm}
 )=\mp k;$
\item[(d)]
 \vspace{1mm}
   \renewcommand{\arraystretch}{2}
  ${\rm tw}(\psfig{figure=zerounoriented.eps,height=0.4cm}
 )-{\rm tw}(\psfig{figure=ktwist.eps,height=0.4cm}
 )=\left\{
   \begin{array}{lr}
   4\lambda -k &
   {\rm in\ the\ case\ (i)} \\
  4\lambda +k &
   {\rm in\ the\ case\ (ii)} \\
   \end{array}
    \right. $
   \par
   \vspace{2mm}
  so the equality (d) holds.
\end{enumerate}
\end{proof}

When one substitutes $a=t^{-3/4}$, $x=-(t^{1/4}+t^{-1/4})$ in the Kauffman polynomial 
 $F(a,x)$ one gets the Jones polynomial $V(t)$ 
(\cite{Li}, see  also \cite{P-1}). Corollary 1.15 gives, 
therefore, some information about the behaviour of $V(t)$ under $t_k$ and $\bar t_k$ moves. 
It happens, however, that one gets no new information comparing with Corollaries1.4, 1.10,  and  Theorem 1.13.

Theorem 1.1 and 1.14 can be stated as one theorem if one uses 
the three variable polynomial $J_L(a,x,z)$ which generalizes 
the Jones-Conway and Kauffman polynomials (see \cite{P-1}), 
however, one cannot gain any new information from this approach.

%\vspace{30pt}
\newpage

\section{Historical background (Fox congruence classes).}
The unknotting number of a knot was considered probably before 
knot theory became a science. It was a natural question to ask how 
many times one has to "cheat" to get from a knot an unknot. 
K.~Reidemeister wrote in 1932 in his book \cite{Re}: "It is very easy 
to define a number of knot invariants so long as one is not concerned 
with giving algorithms for their computation ... One can change each 
knot projection into projection of circle by reversing the overcrossings 
and undercrossings at, say, $k$ double points of the projection. 
The minimal number $u(K)$ of these operations, that is, the minimal 
number of self-piercings, by which a knot is transformed into a circle, 
is a natural measure of knottedness".
 
The first interesting results about unknotting number were found by 
H.~Wendt \cite{We} in 1937. Namely Wendt proved that if $u(K)$ is 
the unknotting number of $K$ and  $e_s$ is the minimal number of 
generators of the group $H_1(M_{K}^{(s)},{\rm Z})$, where 
$M_{K}^{(s)}$ is the cyclic, $s$-fold branched cover of $(S^3, K)$  
then $$e_s\leq u(K)(s-1).$$
$t_k$ and $\bar t_k$ moves 
($\psfig{figure=zerounoriented.eps,height=0.4cm} \rightarrow
\psfig{figure=ktwist.eps,height=0.8cm} $) appear to have been first explicitly 
considered by  S.~Kinoshita in 1957 \cite{Kin-1}. who observed that the 
Wendt inequality is also valid if we allow all $t_{2k}$ and ${\bar{t}}_{2k}$ 
moves, not only $t_2$ moves (see Corollary 2.6(b)). The following year, 
1958, R.~Fox \cite{Fo-1} considered twists of knots and congruence of knots 
modulo $(n,q)$; the notion which is closely related, and in some sense 
more general, than $t_{2k}$ and $\bar t_{2k}$ moves.  Congruence 
modulo $(n,q)$ was  chosen so that the Alexander polynomial (or more 
generally Alexander module) is a good tool to study this.

The same year (1958), S.~Kinoshita \cite{Kin-2} used the Fox twists 
to generalize once more the Wendt inequality (see Corollary 2.6). 
The Fox approach is related to ours so we will present it here 
with some details. We follow the Fox paper \cite{Fo-1} taking into account the corrections made by Kinoshita \cite{Kin-3} and 
Nakanishi and  Suzuki \cite{N-S}. I am grateful to K.~Murasugi and  
H.~Murakami for informing me about the Fox paper and about the Kawauchi 
and Nakanishi conjectures.

Consider the following homeomorphism 
$\tau$ of a 3-disk $D^3=\langle 0,1\rangle\times D^2$: $\tau
(t,z)=(t,e^{2\pi it}z)$. It is the natural extension to $\langle 0,1\rangle\times D^2$ of the Dehn twist on the annulus $\langle 0,1\rangle\times \partial D^2$ (see  Fig. 2.1).

\vspace*{0.8in} \centerline{\psfig{figure=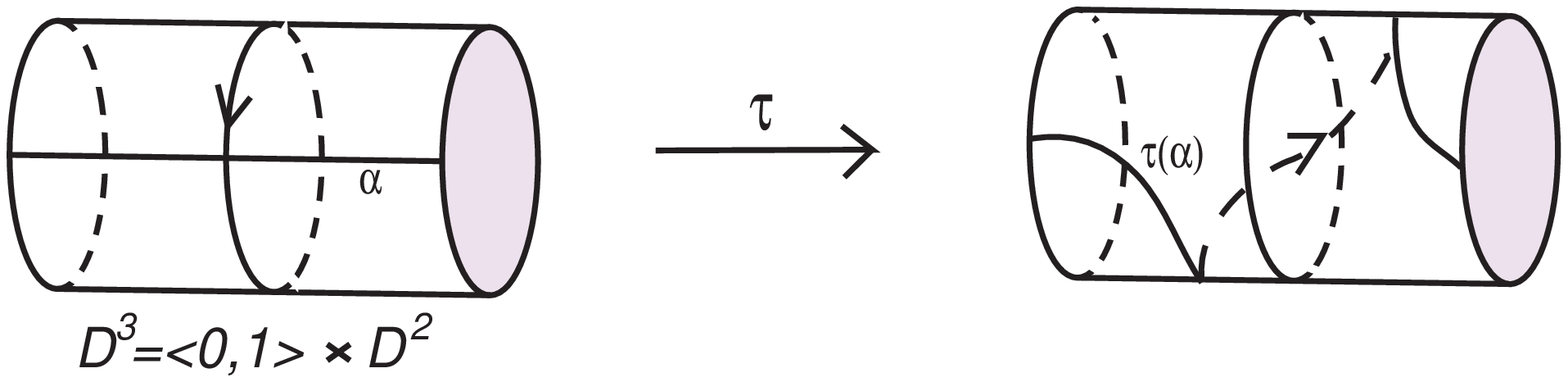,height=4.5cm}}
\begin{center}
Fig. 2.1
\end{center}

We call $\tau$ a simple twist or a Dehn twist. Now whenever we have a properly embedded 2-disk in a 3-manifold $M$ (and either $M$ or a tubular neighbourhood of the disk is oriented), we have uniquely (up to isotopy) associated with the disk the Dehn twist (the twist is carried by a tubular neighbourhood of the disk). In particular for an oriented solid torus there is only one nontrivial Dehn twist, because there is only one, up to isotopy, nontrivial proper disk in it.

Now let $L$ be a link in $\Sigma^3$ (we will assume $\Sigma^3=S^3$, but in fact $\Sigma^3$ can be any homology $3$-sphere), and $D^2$ a disk which cuts $L$ transversely. Let $V_2$ be the solid torus - a small tubular neighbourhood of $\partial D^2$ in $\Sigma^3$, and $V_1$ the closure of its complement 
($V_1=\overline{\Sigma^3-V_2}$). If $\Sigma^3=S^3$, $V_1$
 is a solid torus too. Now perform the Dehn twist on $V_1$ using the disk $D^2$.~The twist restricted to the link $L$ is denoted by $t_{2,q}$ where $q\geq 0$ is the absolute value of the crossing number of $D^2$ and $L$. By $t_{2n,q}$ we denote $t_{2,q}^n$. Notice that our $t_{2n}$ move is special case of $t_{2n,2}$ move, and $\bar
t_{2n}$ move is a special case of $t_{2n,0}$ move. Two oriented links  $L_1$ and $L_2$ are called, by Fox, {\it congruent modulo 
 $n,q$} ($L_1\equiv L_2\ ({\rm mod\ }n,q)$) if one can go from  $L_1$ to $L_2$ using  $t_{2n,q'}^{\mp 1}$, moves (and isotopy), where $q'$ can vary but is always a multiple of $q$. If we allow only $t_{2n,q}^{\mp 1}$ moves then we say, after Nakanishi and Suzuki, that $L_1$ and $L_2$ are $q$-congruent modulo $n$ 
($L_1\equiv_q L_2\ ({\rm mod\ }n)$) or that they are $t_{2n,q}$ equivalent ($L_1\sim_{t_{2n,q}}L_2$).

The Alexander polynomial (and module) is a nice tool for distinguishing nonequivalent links because $L$ and $t_{2n,q}(L)$ are the same outside the ball in which the move occurs.
 \vspace{20pt}

 \begin{theorem}
\begin{enumerate}
\item[(a)]
$t_{2n,q}$ equivalent links have the same Alexander module modulo 
$\frac{(t-1)(t^{nq}-1)}{t^q-1}$, in particular 

\item[(b)]
for $t_{0}^{nq}=1$ ($t_{0}^{q}\neq 1$ or $t_0=1$) 
$\Delta_L(t)\equiv\Delta_{t_{2n,q}(L)}(t)$.
\end{enumerate}
 \end{theorem}

 It can be understood as follows:
 $\Delta_L(t)$ and $\Delta_{t_{2n,q}(L)}(t)$ are equal as elements of the ring 
  $R={\bf Z}[t^{\mp 1}]/(t-1)(1+t^q+\ldots
+t^{(n-1)q})$, up to multiplication by invertible elements of $R$ (in fact up to multiplication by classes of invertible elements in ${\bf
Z}[t^{\mp 1}]$ i.e. $\mp t^p$).

If we substitute $a=i$ and $p=it^{1/2}$ in the Jones-Conway polynomial then we get the (normalized) Alexander polynomial $\Delta (t)\in {\bf Z}[t^{\mp 1}]\cup\sqrt{t}{\bf
Z}[t^{\mp 1}]$. From our Corollary 1.3 follows that if $t^{2k}=1$, 
$t\neq -1$, then $t_{2k}$ move changes  $\Delta_L(t)$ by the factor $\varepsilon =t^k=\mp 1$ (i.e. 
$\Delta_{t_{2k}(L)}(t)=\varepsilon\Delta_L(t)$). Therefore $t_{2k}$ moves have, more less, the same influence on  $\Delta_L(t)$ as more general $t_{2k,2}$ moves; 
however it is not true that every $t_{2k,2}$ move is a combination of $t_{2k}$ moves (see Example 3.8(b)).

\vspace{20pt}
\begin{proof}
Consider a small ball  $B^3$ in which $t_{2n,q}$ move takes place (Fig. 2.2). 
$B^3\cap L$ consists of $m$ parallel strings.

\vspace*{0.8in} \centerline{\psfig{figure=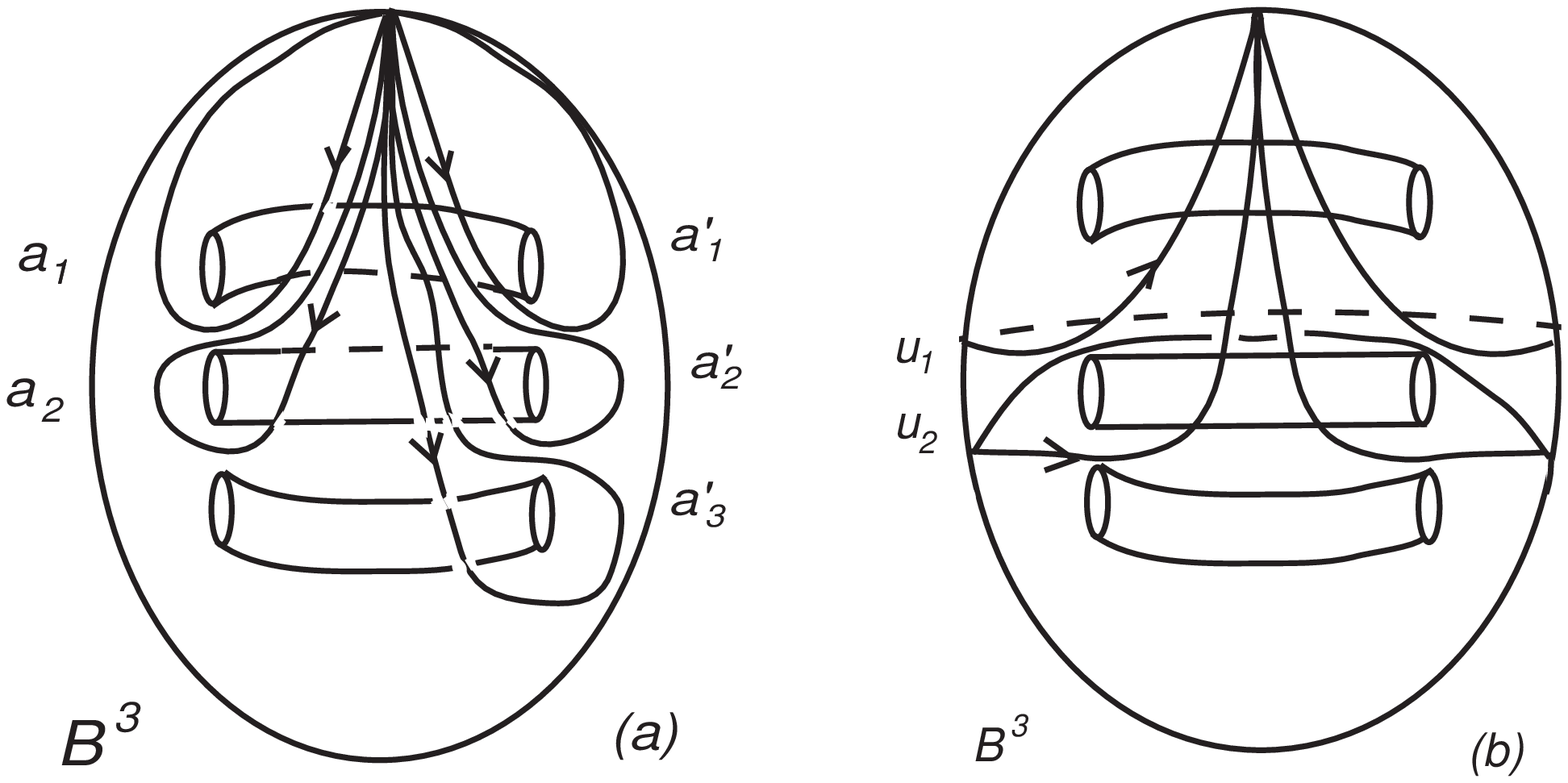,height=4.5cm}}
\begin{center}
$m=3$
\par
Fig. 2.2
\end{center}

$\Sigma^3-L-{\rm int}B^3$ is homeomorphic to $\Sigma^3-t_{2n,q}(L)-{\rm int}B^3$~ and the fundamental groups of these spaces have the following presentation:

$$\{a_1,\ldots ,a_{m-1}, a'_1,\ldots ,a'_{m}, x_1,\ldots
:r_1,\ldots\},$$

where $a_1,\ldots ,a_{m-1}, a'_1,\ldots ,a'_{m}$ form a basis of the free group  $\pi_{1}(\partial B^3-L)$; see Fig. 
2.2(a).

$\Sigma^3-L$ and $\Sigma^3-t_{2n,q}(L)$ can be obtained from $\Sigma^3-L-B^3$ by adding $m-1$ two-disks in the appropriate way (Fig. 2.2(b)). Therefore 
$$\pi_{1}(\Sigma^3-L)=\{a_1,\ldots ,a_{m-1}, a'_1,\ldots ,a'_{m},
x_1,\ldots : u_1, u_2,\ldots , u_{m-1}, r_1, \ldots\},$$

where $u_i=a'_ia_{i}^{-1}$, $i=1,\ldots ,m-1$ (see Fig. 2.2(b)),
 and

$$\pi_{1}(\Sigma^3-t_{2n,q}(L))=\{a_1,\ldots ,a_{m-1}, a'_1,\ldots
,a'_{m}, x_1,\ldots : \tau (u_1), \tau (u_2),\ldots , \tau (
u_{m-1}), r_1, \ldots\},$$

where  $\tau (u_i)=(a'_m,\ldots ,a'_1)^na'_i(a'_m,\ldots
,a'_1)^{-n}a_{i}^{'-1}$.

Consider the natural projections 
$p=p_2p_1:\pi_{1}(\Sigma^3-L)\stackrel{P_1}\rightarrow  H_1
(\Sigma^3-L)\stackrel{P_1}\rightarrow {\bf Z},$ where  $p$ 
sends meridians of $L$ onto $t$ - a generator of integers, and 
$p':\pi_{1}(\Sigma^3-t_{2n,q}(L))\rightarrow
{\bf Z}$. Then 

$p(a_i)=p(a'_i)=t^{\mp 1}$ and $p(a'_m a'_{m-1}\ldots
a'_1)=t^q=p(a_m a_{m-1}\ldots a_1)$ (without the loss of generality one can assume that the crossing number of $D^2$ and  $L$ is nonnegative so equal to $q$).
\par
In particular if $i$ and $i'$ are embeddings of $\Sigma^3-L-B^3$, in $\Sigma^3-L$ and 
$\Sigma^3-t_{2n,q}(L)$ respectively then $pi_{\ast}=p'i'_{\ast}$ (lack of this condition was the source of the mistake in the Fox paper \cite{Fo-1}).

Now one can use Fox calculus to find Alexander-Fox modules of group representations 
$p:\pi_1(\Sigma^3-L)\rightarrow {\bf Z}$ and
$p':\pi_1(\Sigma^3-t_{2n,q}(L))\rightarrow {\bf Z}$, and because

 \vspace{1mm}
   \renewcommand{\arraystretch}{2}
   \begin{center}
  $p_{\ast}(\frac{\partial u_i}{\partial a'_j})
=\left\{
   \begin{array}{lr}
   0 &  {\rm if }~~~  i\neq j \\
 1 & {\rm if }  ~~~i=j
   \end{array}
    \right. $ and
    \end{center}
   \par
   \vspace{2mm}

 \vspace{1mm}
   \renewcommand{\arraystretch}{2}
  $p'_{\ast}(\frac{\partial\tau (u_i)}{\partial a'_j})
=\left\{
   \begin{array}{lr}
   \frac{(1-t)(1-t^{nq})t^{b_j}}{1-t^q},
   {\rm\  where~ } ~ t^{b_j}=p(a'_m\ldots a'_{j+n}) & {\rm if}~~~ i\neq j \\
 \frac{(1-t)(1-t)^{nq}t^{b_j}}{1-t^q}+t^{nq} &
   {\rm if}~~~  i=j
   \end{array}
    \right. $
   \par
   \vspace{2mm}

therefore

$p_{\ast}(\frac{\partial u_i}{\partial
a'_j})=p'_{\ast}(\frac{\partial\tau (u_i)}{\partial a'_j}){\rm\
mod\ }\frac{(1-t)(1-t^{nq})}{1-t^q}$ and one gets:

\vspace{20pt}

\begin{lemma}
The Alexander-Fox modules of 
$p:\pi_1(\Sigma^3-L)\rightarrow {\bf Z}$ and 
$p':\pi_1(\Sigma^3-t_{2n,q}(L))\rightarrow {\bf Z}$ 
can be represented by the matrices which are the same modulo 
$\frac{(1-t)(1-t^{nq})}{1-t^q}$.
\end{lemma}

\vspace{20pt}

Theorem 2.1 follows immediately from the lemma.

\end{proof}

\vspace{20pt}

\begin{corollary}
\begin{enumerate}
\item[(a)]
Let  $qn$ be a multiple of  $s$ and 
$k=\frac{n\cdot {\rm gcd}(q,s)}{s}$; where gcd() is the greatest common divisor, then a $t_{2n,q}$ move does not change $H_1(M_{L}^{(s)},{\bf Z}_k)$. 
In particular 

\item[(b)]
If  $q$ is a multiple of  $s$ (e.g. $q=0$) then a 
$t_{2n,q}$ move does not change $H_1(M_{L}^{(s)},{\bf Z}_n)$.

\item[(c)]
Let $n$ be a multiple of $s$, and $s$ and  $q$ are coprime then a $t_{2n,q}$ move does not change  $H_1(M_{L}^{(s)},{\bf Z})$.
\end{enumerate}
\end{corollary}

\begin{proof}
Alexander matrices can be used to describe $H_1(M_{L}^{(s)},{\bf Z})$ as ${\bf Z}[{\bf Z}_s]$ module 
 ($t^s=1$). Then we use  Lemma 2.2.
\end{proof}

\vspace{20pt}

\begin{corollary}
\begin{enumerate}
\item[(a)]
Let $\overline{\overline u}_m(L)$ denote the minimal number of $t_{2n,q}$ moves (we allow different $n$ or $q$) but the number of strings involved in a $t_{2n,q}$ must be less or equal $m$) which are needed to change a given link  $L$ into unlink then 

$$|e_s-s(c(L)-1)|\leq (s-1)(m-1)\overline{\overline u}_m(L),$$

where  $c(L)$ is the number of components of $L$ and  $e_s$ is the minimal number of generators of 
$H_1(M_{L}^{(s)},{\bf Z})$. In particular for $\bar
u(L)=\overline{\overline u}_2(L)$ one gets:

\item[(b)]
\cite{Kin-1} The minimal number of $t_{2n}$ or $\bar t_{2n}$ moves which are needed  to change a given link $L$ into unlink ($\bar u(L)$), satisfies:

$$|e_s-s(c(L)-1)|\leq (s-1)\bar u(L).$$

\end{enumerate}
\end{corollary}

\begin{proof}
If  $T_n$ is a trivial link of $n$ components then 
$e_s(T_n)=s(n-1)$. By the proof of Lemma 2.2, 
$H_1(M_{t_{2n,q}L}^{(s)},{\bf Z})$ has a presentation which differs from a presentation of $H_1(M_{L}^{(s)},{\bf Z})$ at most in $(s-1)(m-1)$ rows (we use additionally the fact that 
$\frac{t^s-1}{t-1}=1+\ldots t^{s-1}$ is an annihilator of
$H_1(M_{L}^{(s)},{\bf Z})$ (see \cite{B-Z}), so 

$$|e_s(t_{2n,q}(L))-e_s(L)|\leq (s-1)(m-1).$$
\end{proof}

\vspace{20pt}

The Fox method  (and Lemma 2.2) can be modified so that one can get the result about $t_k$ moves, for $k$ odd, analogous to Corollary 2.3 (compare \cite{Ki}).

Consider a small ball $B^3$ in which a $t_k$ move takes place (Fig. 2.3).

\vspace*{0.8in} \centerline{\psfig{figure=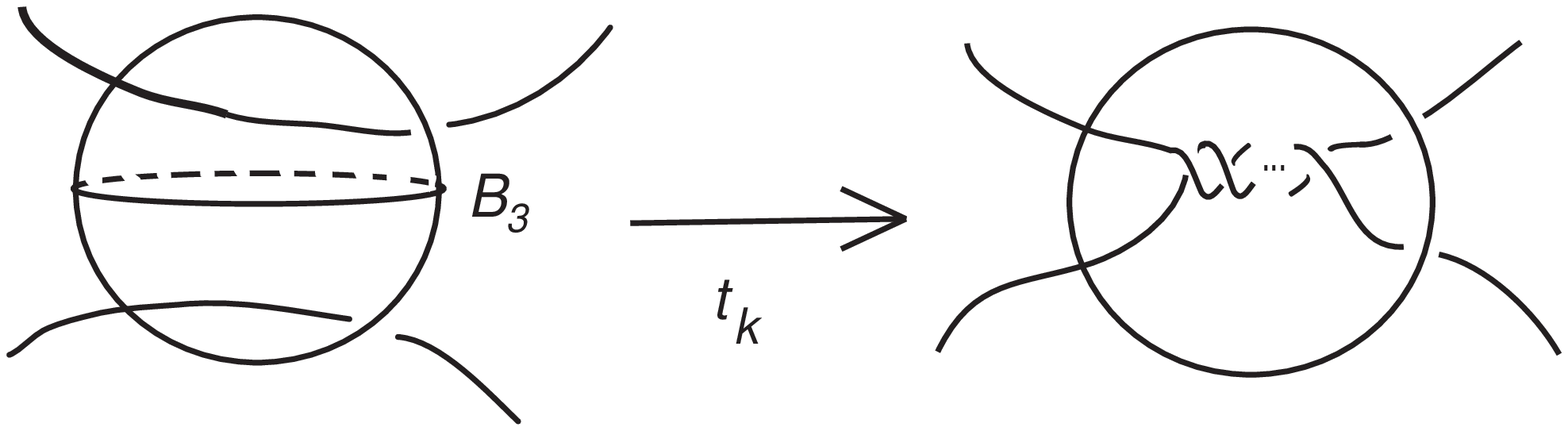,height=4.2cm}}
\begin{center}
Fig. 2.3
\end{center}

$\Sigma^3-L-{\rm int}B^3$ is homeomorphic to 
$\Sigma^3-t_k(L)-{\rm int}B^3$ and the fundamental groups of them have the following presentation:

$$\{a,b,c,x_1,x_2,\ldots :r_1, r_2\},$$

where  $a$, $b$ and $c$ are classes of curves (generators of $\pi_1(\partial B^3-L)$) shown on Fig. 2.4.

%\newpage

\vspace*{0.4in} \centerline{\psfig{figure=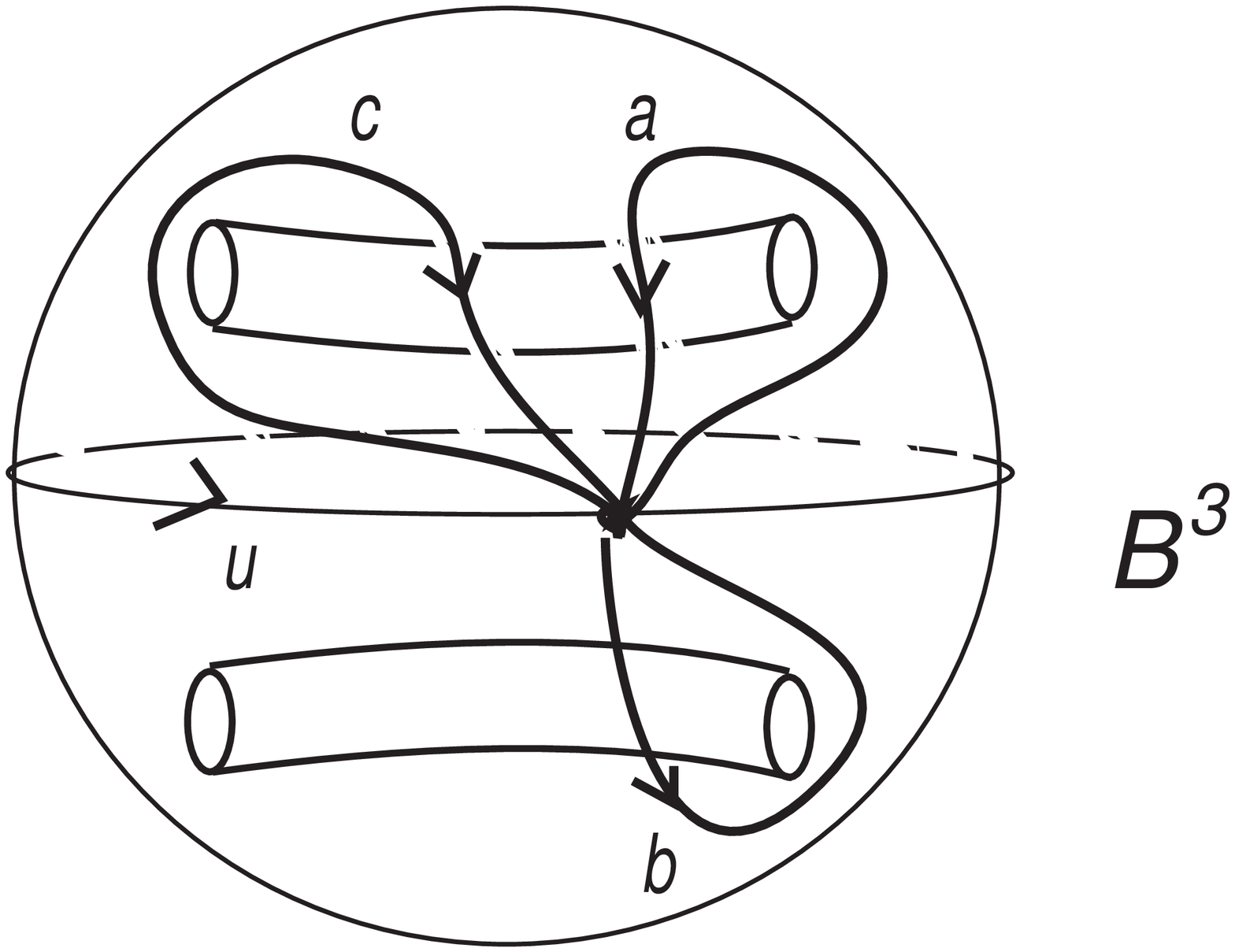,height=4.3cm}}
\begin{center}
Fig. 2.4
\end{center}

\newpage

If we add a $2$-handle along $u=ac^{-1}$ we get $\Sigma^3-L$ 
and if we add a 2-handle along $(ba)^{(k-1)/2}b(ba)^{-(k-1)/2}c^{-1}$ we get 
$\Sigma^3-t_k(L)$; therefore 
$$\pi_1(\Sigma^3-L)=\{a,b,c,x_1,x_2,\ldots :ac^{-1}=1,r_1,
r_2,\ldots \}, $$

$$\pi_1(\Sigma^3-t_k(L))=\{a,b,c,x_1,x_2,\ldots
:(ba)^{(k-1)/2}b(ba)^{-(k-1)/2}c^{-1},r_1, r_2,\ldots \}. $$

Consider the natural projections 
$p:\pi_1(\Sigma^3-L)\rightarrow{\bf Z}$ and
$p':\pi_1(\Sigma^3-t_n(L))\rightarrow{\bf Z}.$

Then we have:

$$p(a)=p(b)=p(c)=t,\ p'(a)=p'(b)=p'(c)=t.$$

Now we calculate that ($r_0=(ba)^{(k-1)/2}b(ba)^{-(k-1)/2}c^{-1}$):

$$p(\frac{\partial ac^{-1}}{\partial a})=1,\ p'(\frac{\partial
r_0}{\partial a})=1-\frac{t^k+1}{t+1},$$

$$p(\frac{\partial ac^{-1}}{\partial b})=0,\ p'(\frac{\partial
r_0}{\partial b})=\frac{t^k+1}{t+1},$$

$$p(\frac{\partial ac^{-1}}{\partial c})=-1,\ p'(\frac{\partial
r_0}{\partial c})=-1,$$

and we get:

\vspace{20pt}

\begin{lemma}
The Alexander-Fox modules of $L$ and $t_k(L)$ can be presented by the following matrices.

\begin{tabular}{c|cccccc}
$L$ & $a$ & $b$ & $c$ & $x_1$ & $x_2$ & $\ldots $\\

 \hline

$ac^{-1}=1$ & $1$ & $0$ & $-1$ & $0$ & $0$ & $\ldots$\\

$r_1 $ & $\ast$ & $\ast$ & &  $\ast$ & $\ast $ & \\

$\ldots$ & & &  $\ldots$ & & &

\end{tabular}\ \ \
\begin{tabular}{c|cccccc}
$t_k(L)$ & $a$ & $b$ & $c$ & $ x_1$ & $x_2$ & $\ldots $\\

 \hline

$r_0$ & $1-\frac{t^k+1}{t+1}$ & $\frac{t^k+1}{t+1}$ & -1 & 0 & 0 &
$\ldots$\\

$r_1 $ & $\ast$ & $\ast$ & &  $\ast$ & $\ast $ & $\ast $ \\

$\ldots$ & & &  $\ldots$ & & &

\end{tabular}
\end{lemma}

\vspace{20pt}

\begin{corollary}
\begin{enumerate}
\item[(a)]
For  $\frac{t^k+1}{t+1}=0$ ($k$-odd)
$t_k$-equivalent links have the same Alexander module, in particular

\item[(b)]
$\Delta_{t_k(L)}(t)\equiv \mp t^i\Delta_L(t) ({\rm \ mod\
}\frac{t^k+1}{t+1}).$

In fact from Corollary 1.3 follows that for a normalized Alexander polynomial 
$\Delta_L(t)\equiv \mp
i\Delta_{t_k(L)}(t)({\rm \ mod\ }\frac{t^k+1}{t+1})$ or precisely 

$$\Delta_L(t)\equiv t^{k/2}\Delta_{t_k(L)}(t)({\rm \ mod\
}\frac{t^k+1}{t+1}).$$
\end{enumerate}
\end{corollary}

 We can slightly generalize the results of Wendt and Kinoshita using Lemma 2.5.

\vspace{20pt}

\begin{corollary}
Let $\overline{\overline{\overline u}}_n(L)$ denote the minimal number of $\bar{t_{2k}}$ or $t_k$ moves which are needed to change a given oriented link $L$ into unlink of $n$ components, then 

$$|e_s-s(n-1)|\leq (s-1)\overline{\overline{\overline u}}_n(L),$$

where  $e_s$ is the minimal number of generators of $H_1(M_{L}^{(s)},{\bf Z})$.
\end{corollary}

\vspace{30pt}

\section{Applications and Speculations}

We start this part by proving two "folklore" results which link Goeritz and Seifert matrices with $t_k$ or ${\bar{t}}_k$ moves.

\vspace{20pt}

\begin{theorem}
\begin{enumerate}
\item[(a)]
There exist Goeritz's matrices for $L$ and $t_k(L)$ (or $\bar t_k(L)$) which are the same modulo $k$.
\item[(b)]
$t_k$ and  $\bar t_k$-moves preserves $H_1(M_{L}^{(2)},{\bf
Z}_k)$.
\end{enumerate}
\end{theorem}

\begin{proof}
For the convenience we start from the definition of Goeritz's matrix (\cite{Goe,Gor}). Colour the regions of the diagram of an unoriented link alternately black and white, the unbounded region  $X_0$ being coloured white, and number the other white regions $X_1,\ldots ,X_n$. Assign an incidence number $\eta (p)=\mp 1$  to each crossing point $p$ as shown in Fig. 3.1. Then define $n\times n$ Goeritz's  matrix $G=(g_{ij})$ by 

\vspace*{0.8in} \centerline{\psfig{figure=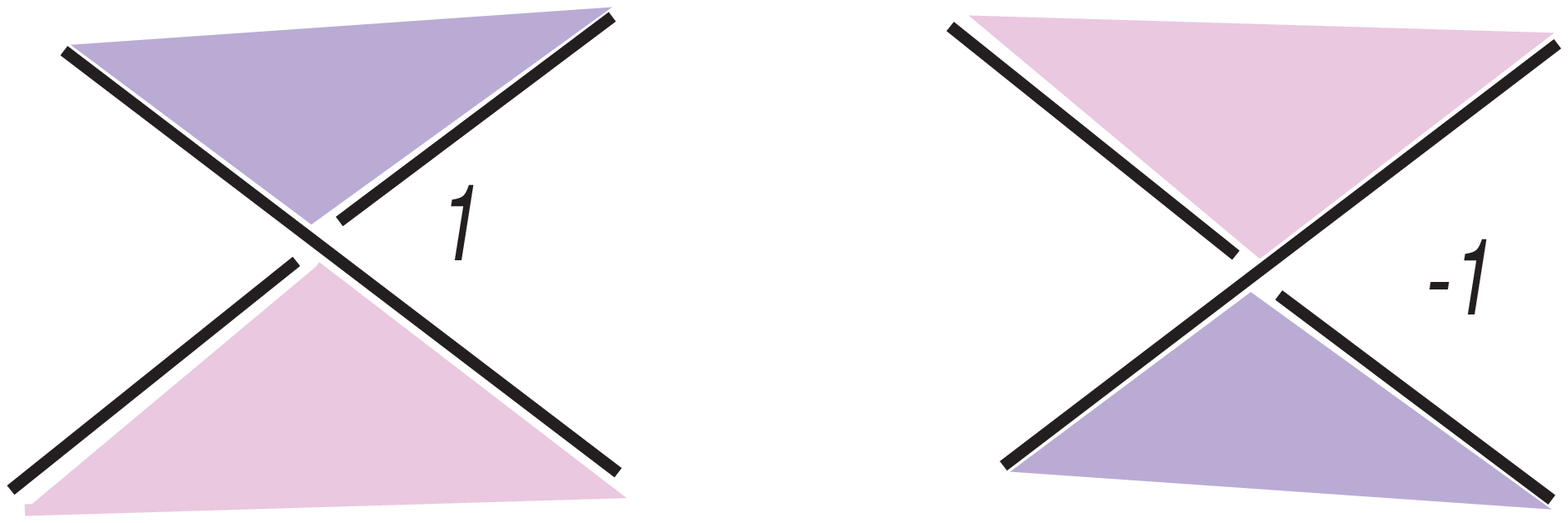,height=3cm}}
\begin{center}
Fig. 3.1
\end{center}

\vspace{20pt}
   \renewcommand{\arraystretch}{1}
   $g_{ij}=\left\{
  \begin{array}{lr}
   \sum\eta (c)  {\rm\ \ \ \  summed~ over~ crossings~ points~ p~ adjacent~ to~~~} X_{\rm\ i\ } ~{\rm\ and}  ~X_{\rm\ j\ }~\\
   ~~~~~~~~~~{\rm{if}} ~i\neq j\ (i\cdot j \geq 1) \\
   -\sum\eta (c)
   {\rm\ \ \ \  summed~ over~ crossings~ points~} p~ {\rm{adjacent}}~ {\rm{to}}~
   \\
    ~~~~~~~~~~X_i~ {\rm{and~ to~ some}}~ X_j ~(i\neq j)~ {\rm{if}}~ i=j ~{\rm{(}}i \geq 1{\rm{)}}. \\
   \end{array}
   \right.$
   \par
   \vspace{20pt}

  Now consider the Fig. 3.2 with white regions $X_i$ and $X_j$.

\vspace*{0.8in} \centerline{\psfig{figure=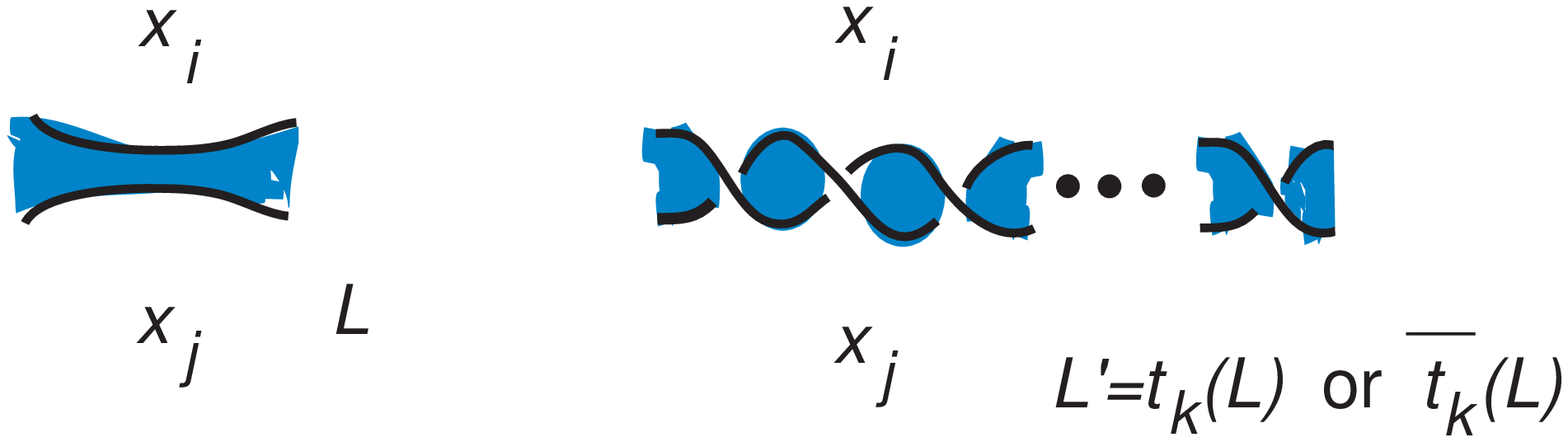,height=4.5cm}}
\begin{center}
Fig. 3.2
\end{center}

There are two possible cases:

\begin{enumerate}
\item[(i)]
$X_i=X_j$, then  $G_L=G_{L'}$; in fact $L$ is isotopic to $L'$.

\item[(ii)]
$X_i\neq X_j$, we can assume that $i=0$ and  $j=1$ then

\vspace{20pt}

   $$G_{L'}=\left[
  \begin{array}{cccc}
  g_{11}+k, & g_{12} & \ldots , & g_{1n} \\
   & \cdot\cdot\cdot & &  \\
   g_{n1}, & g_{n2} & \ldots , & g_{nn}
   \end{array}
   \right]~~where~~ G_L=(g_{ij})$$.
   \par
\end{enumerate}

The part (b) of the Theorem 3.1 follows from the fact that $G_L$ is a presentation matrix for $H_1(M_{L}^{(2)},{\bf Z})$.
\end{proof}

\vspace{20pt}

An alternative proof of (b) can be given by considering Dehn surgery on 
$M_{L}^{(2)}$ corresponding to  $t_k$ or $\bar t_k$ move on $L$.

\vspace{20pt}

\begin{theorem}
\begin{enumerate}
\item[(a)]
Consider a  $t_{2k,0}$ move of Fox (e.g. $\bar t_{2k}$ move), then there exist Seifert matrices for $L$ 
 and  $t_{2k,0}(L)$ which are the same modulo $k$.

 \item[(b)]

$t_{2k,0}$ move preserves $H_1(M_{L}^{(s)},{\bf Z}_k)$ for any $s$
.
\end{enumerate}
\end{theorem}

\begin{proof}
One can find a Seifert surface $S$ for $L$ which cuts the disk $D^2$ which supports the 
$t_{2k,0}$ move, as shown in Fig. 3.3. Then the Seifert matrix for $L$ defined by $S$ and for 
$t_{2k,0}(L)$ defined by $t_{2k,0}(S)$ satisfy the condition (a).

 \vspace*{0.8in}
\centerline{\psfig{figure=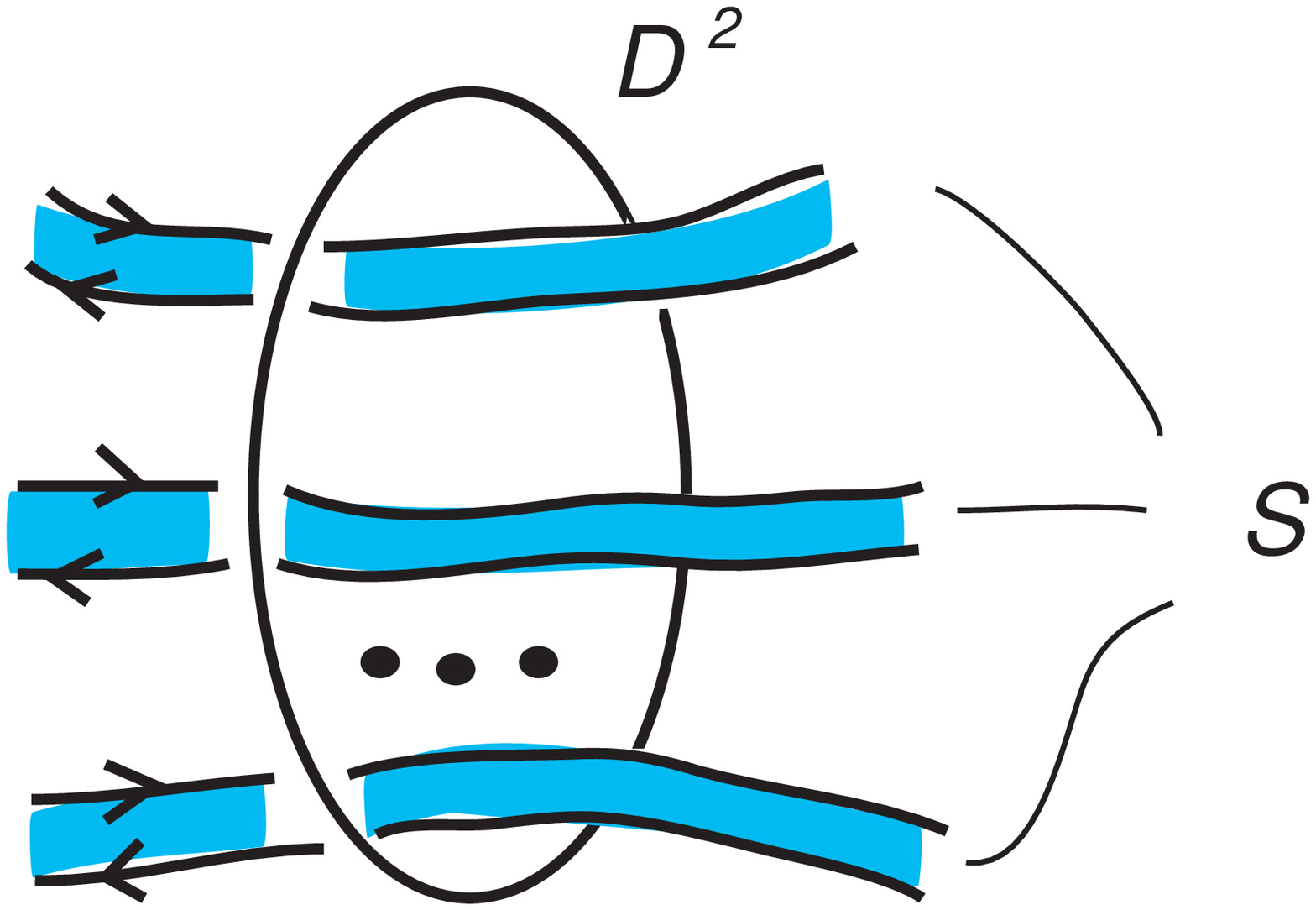,height=4.5cm}}
\begin{center}
Fig. 3.3
\end{center}

(b) follows from (a) because a presentation matrix for 
$H_1(M_{L}^{(s)},{\bf Z})$ can be built of blocks of the shape $\mp V$, 
$\mp V^T$, $\mp (V+V^T)$, where $V$ is a Seifert matrix of $L$ and $V^T$ its transpose. 
On the other hand, (b) is a special case of Theorem 2.3(b).
\end{proof}

\vspace{20pt}

\begin{example}
\begin{enumerate}
\item[(a)]
The trivial knot $(T_1)$ and the (right handed) trefoil knot 
$(3_1)$ are $t_4$ equivalent. The figure eight knot  $(4_1)$ and the $5_2$ knot are 
$t_4$ equivalent however they are not $t_4$ equivalent to $T_1$ or $3_1$.

\item[(b)]
$T_1$ and  $5_2$ are $\bar t_4$ equivalent. $3_1$ and $4_1$ 
are $\bar t_4$ equivalent but they are not $\bar
t_4$ equivalent to  $T_1$ or $5_2$.
\end{enumerate}

First parts of  (a) and  (b) are illustrated in Fig. 3.4.

\vspace*{0.8in} \centerline{\psfig{figure=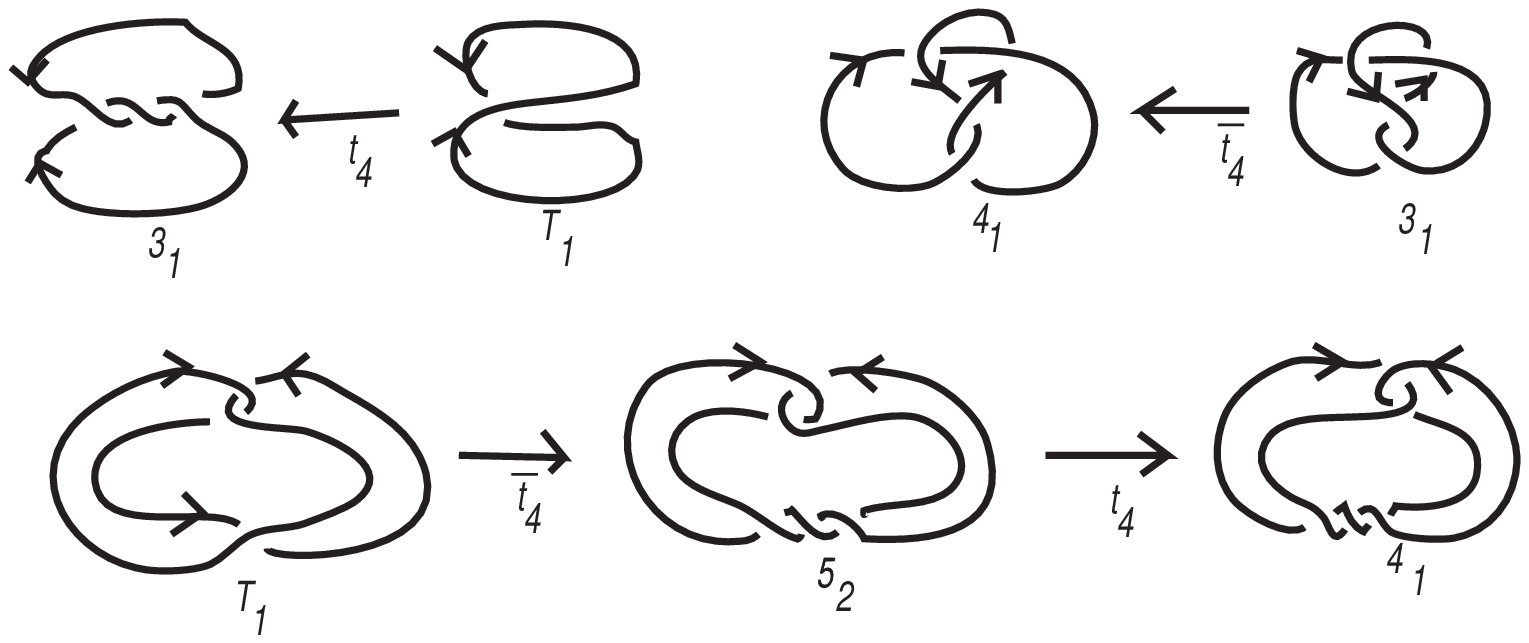,height=5.5cm}}
\begin{center}
Fig. 3.4
\end{center}

The second parts follow from Corollaries 1.2 and 1.8 ($t_4$ move changes $P_L(a, \sqrt{2})$  by the factor 
$-a^{-4}$ and $\bar t_4$ preserves  $P_L(1,z)$) and the following computation:

$$P_{T_1}(a,z)=1$$

$$P_{3_1}(a,z)=-a^{-4}-2a^{-2}+z^2a^{-2}; \ P_{3_1}(a,
\sqrt{2})=-a^{-4},\ P_{3_1}(1,z)=z^2-3$$

$$P_{4_1}(a,z)=-a^{-2}-1-a^2+z^2;\ P_{4_1}(a,
\sqrt{2})=-a^{-2}+1-a^2;\ P_{4_1}(1,z)=z^2-3$$

$$P_{5_2}(a,z)=-a^{-2}+a^4+a^6+z^2(a^2-a^4);\ P_{5_2}(a,
\sqrt{2})=-a^4(-a^{-2}+1-a^2);\ P_{5_2}(1,z)=1.$$

\end{example}

\vspace{20pt}

\begin{example}
Every closed 3-braid knot is $t_4$ equivalent to the trivial knot or the figure eight knot.
It is not an unexpected result because the quotient group $B_3/(\delta_{1}^{4})$ is finite \cite{Cox}. 
In fact a calculation shows that $B_3/(\delta_{1}^{4})$ has only two classes of knots 
(represented by $T_1$ and  $4_1$). 
Because all presentations of  $4_1$  as a 3-braid (e.g. $\delta_1\delta_{2}^{-1}\delta_1\delta_{2}^{-1}$) 
have the same exponent sum (equal to 0) therefore for every knot $K$ which is $t_4$ equivalent to $4_1$,
 each of its presentation as a 3-braid has the same exponent sum (equal to  $4|4_1,K|_{t_4}^{\rm lev}$; compare Corollary 1.2). 
More in this direction can be got using other $t_k$ moves, compare Example 3.11, however it has been generally proved by  H.Morton \cite{Mo} and J.Birman that if $L$ is not a $(2,k)$ torus link then the exponent sum of $L$ does not depend on the presentation of $L$ as a closed $3$-braid.
\end{example}

\vspace{20pt}

\begin{example}
Consider the following theorem of H.Murakami \cite{Mur-1} (see also \cite{L-M-2}):

\vspace{1mm}
   \renewcommand{\arraystretch}{2}
   $P_L(1,\sqrt{2})=V_L(i)=\left\{
   \begin{array}{lr}
 (\sqrt{2})^{c(L)-1}(-1)^{{\rm Arf}(L)} &
   {\rm if\ }{\rm Arf}(L){\rm\ exists}\\
   0 &
   {\rm otherwise\ },\ \\

   \end{array}
   \right. $
   \par
   \vspace{2mm}

where $c(L)$ denotes the number of components of $L$, ${\rm
Arf}(L)$ is the Arf (or  Robertello) invariant (see \cite{Rob} or \cite{Ka-2}), and $t=i$ in $V_L(t)$ 
should be understood as $t^{1/2}=-e^{\pi i/4}$. Notice that our convention differs slightly 
from that of \cite{L-M-1} or  \cite{L-M-2} namely  $P_L(a,z)=P_L({\ell},
-m)=(-1)^{c(L)-1}P_L({\ell}, m)$.

It follows from Corollaries 1.2 and  1.8 that $t_4$ move changes 
$P_L(1,\sqrt{2})$ by factor  $-1$ and $\bar t_4$ move preserves 
 $P_L(1,\sqrt{2})$. Furthermore for $T_n$ - the 
trivial link of  $n$ components $P_{T_n}(1,\sqrt{2})=(\sqrt{2})^{n-1}$. 
On the other hand the Arf invariant of a trivial link is equal to zero, 
$t_4$ move changes the Arf invariant (if defined) and $\bar t_4$ move
 preserves it (see 
\cite{Ka-2}). Therefore the Murakami theorem follows immediately from 
the above observations for a link which is  $t_4$, $\bar t_4$  equivalent to a trivial link (
i.e. a link which can be obtained from a trivial one using $t_4$ and  $\bar t_4$ moves). 
This should be confronted with the following conjecture:

%\vspace{25pt}
\newpage

\begin{conjecture}(Kawauchi - Nakanishi)
\begin{enumerate}
\item[(a)]
If two links $L_1$ and  $L_2$ are homotopic then they 
are $t_4$, $\bar t_4$ equivalent\footnote{Added for e-print: \ 
The conjecture has been disproved in \cite{D-P-2} for links of three 
or more components. For two component links it is still an open problem 
whether any such link is $t_4$, $\bar t_4$ equivalent to $T_2$ or 
the Hopf link.}. In particular:

\item[(b)]
Every knot is $t_4$, $\bar t_4$ equivalent to the unknot. 
\end{enumerate}
\end{conjecture}

Conjecture 3.6 has been verified for the 2-bridge links, 
closed $3$-string braids and  pretzel links.
\end{example}

\vspace{20pt}

\begin{example}
Consider the following  $t_{\Delta^2}$-move ($\Delta^2$-twist) 
on oriented diagrams of links (Fig. 3.5).

\vspace*{0.8in} \centerline{\psfig{figure=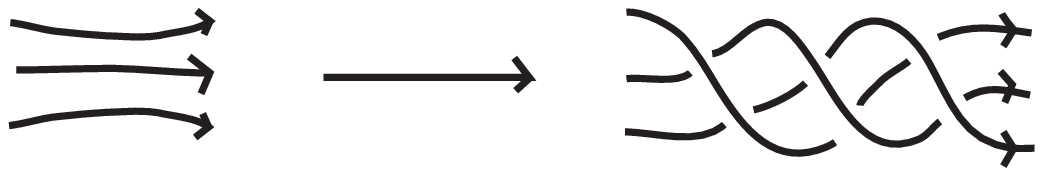,height=2.1cm}}
\begin{center}
Fig. 3.5
\end{center}

A $t_4$ move can be obtained from a $t_{\Delta^2}$-move (and isotopy)
 as it is illustrated in Fig.3.6.

\vspace*{0.8in} \centerline{\psfig{figure=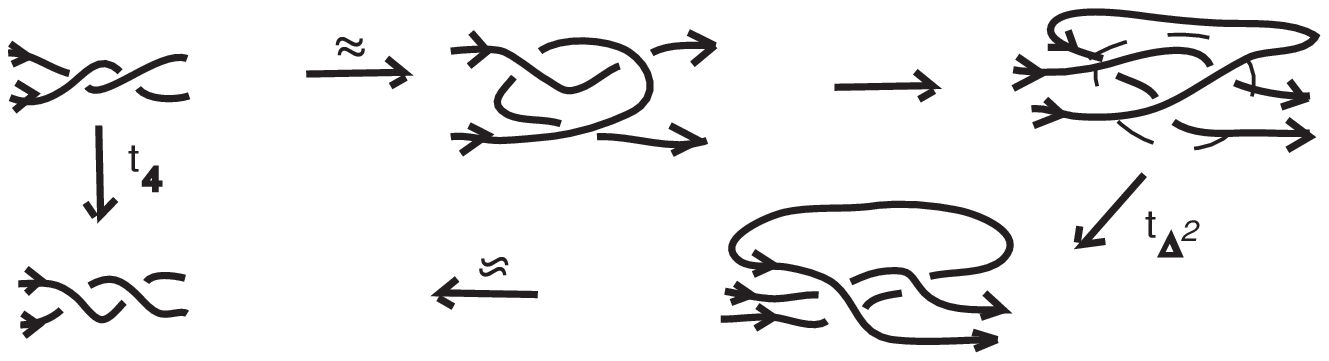,height=4.5cm}}
\begin{center}
Fig. 3.6
\end{center}

J.Birman  and  B.Wajnryb \cite{B-W} have proven that two links are 
$t_{\Delta^2}$ equivalent iff they have the same number of components and the same number of components with odd linking number with the rest of the link. Because  a  $\bar t_4$ move preserves the number of components and all 
linking numbers modulo 2, 
therefore it can be obtained as a combination of 
$t_{\Delta^2}^{\mp 1}$ moves. In fact it follows from \cite{B-W} that in order to get 
$\bar t_4$ move we can always use an even number of 
$t_{\Delta^2}^{\mp 1}$ moves. Furthermore a 
$t_{\Delta^2}$ move changes the Arf invariant 
(if it exists) and therefore 
$V_L(i)=-V{t_{{\Delta^2}(L)}}(i)$. The last equality can be 
also proven elementary without using \cite{B-W}. 
Finally observe that not every  $t_{\Delta^2}$ move is a combination of 
 $t_4$, $\bar t_4$ moves. The reason is that 
$t_4$ and $\bar t_4$ moves preserve all linking numbers mod $2$ 
but it is not always the case for a 
$t_{\Delta^2}$ move (see Fig. 3.7 for an example of links which are 
$t_{\Delta^2}$ equivalent but not $t_4$, $\bar
t_4$ equivalent).

\vspace*{0.8in} \centerline{\psfig{figure=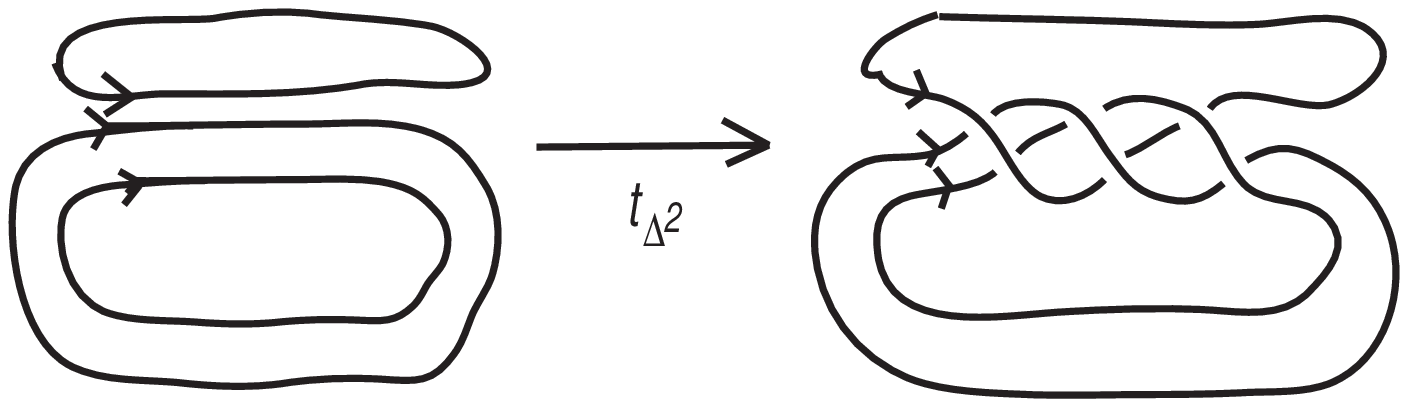,height=4.5cm}}
\begin{center}
Fig. 3.7
\end{center}

\end{example}

\vspace{20pt}

\begin{example}
\begin{enumerate}
\item[(a)]
a $t_{\Delta^2}$ move is a special case of $t_{2,3}$ moves of Fox
but it follows from \cite{B-W} that any $t_{2,3}$ move is a combination of 
$t_{\Delta^2}$ moves. In fact, every $t_{2,3}$ move preserves the number 
of components and the number of components with odd linking number with the rest of the link.
Similarly any $t_{2,2q+1}$ is a combination of $t_{\Delta^2}$ moves. 
\item[(b)] a $t_4$ move is a special case of $t_{4,2}$ moves of Fox.
 There are $t_{4,2}$ equivalent links which are not $t_4$ equivalent. (Fig. 3.8).

\vspace*{0.8in} \centerline{\psfig{figure=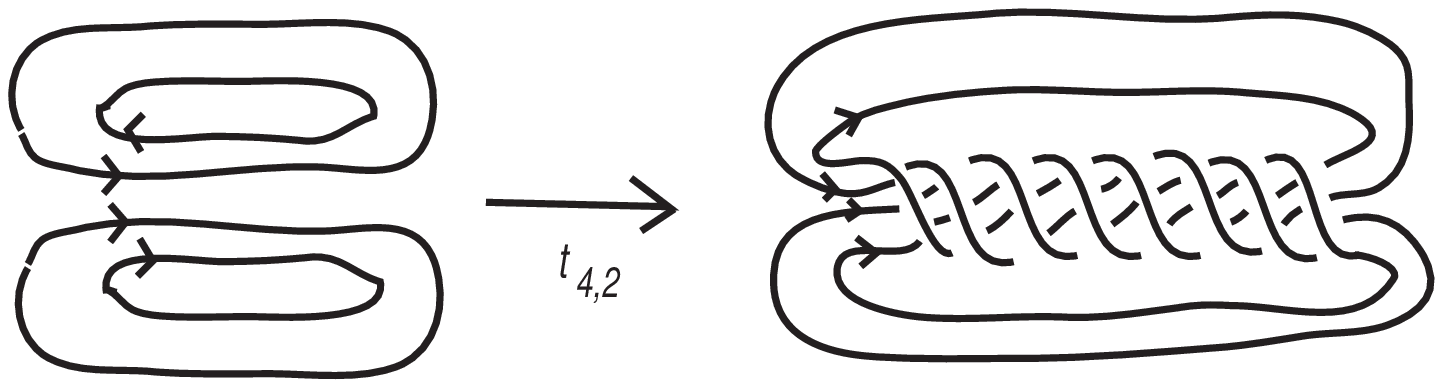,height=4.5cm}}
\begin{center}
Fig. 3.8
\end{center}

Two links of Fig. 3.8 are not $t_4$ equivalent because their sublinks of Fig. 3.9
 are not $t_4$ equivalent.

\vspace*{0.8in} \centerline{\psfig{figure=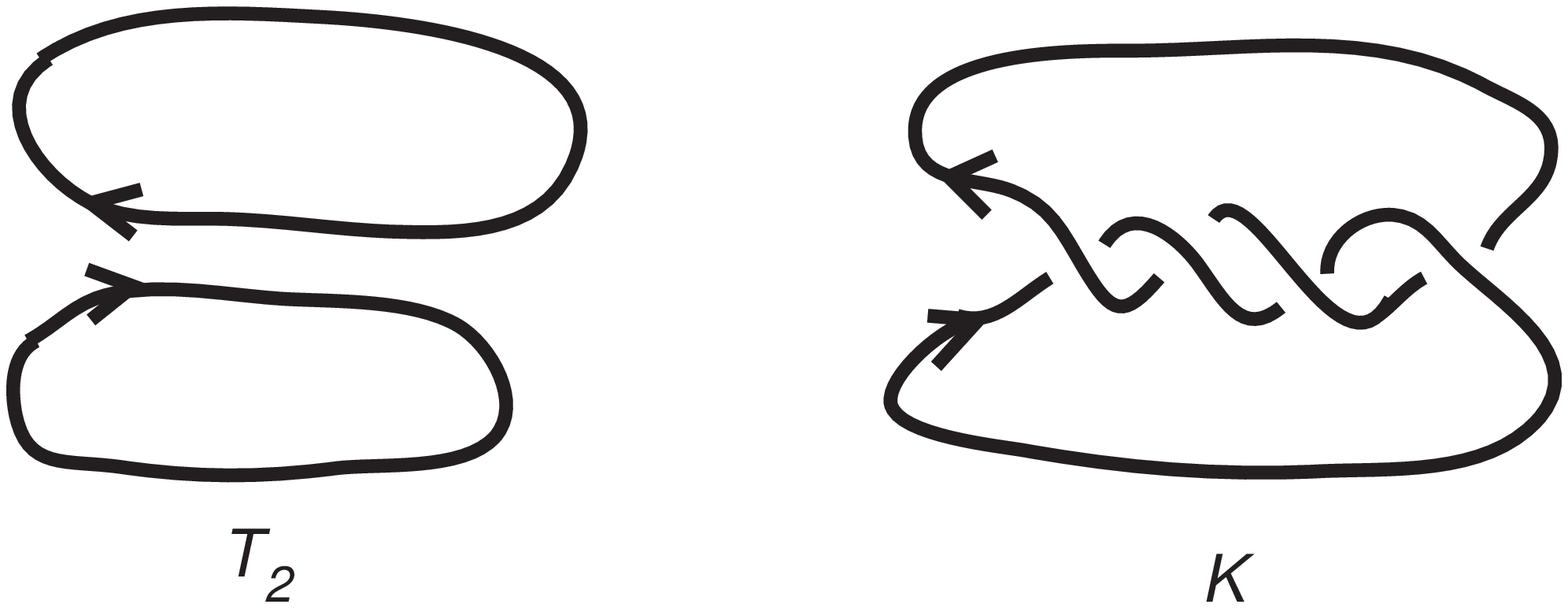,height=4.5cm}}
\begin{center}
Fig. 3.9
\end{center}

Namely  $P_{T_2}(a,\sqrt{2})=\frac{a+a^{-1}}{\sqrt{2}}$ and 
$P_K(a,\sqrt{2})=\frac{a^{-5}-a^{-3}+2a^{-1}}{\sqrt{2}}$, therefore by Corollary 1.2 
 $T_2$ and  $K$  are not $t_4$ equivalent.
\end{enumerate}
\end{example}

\vspace{20pt}

\begin{example}
\begin{enumerate}
\item[(a)]
The square knot ($3_1\#\bar 3_1$), the (right-handed) granny knot 
($3_1\# 3_1$), and $T_3$ (the trivial $3$-component link) are 
$t_3$ equivalent.

\item[(b)]
The trefoil knot ($3_1$) and $T_2$ (the trivial link of~ $2$-components) 
are  $t_3$ equivalent.

\item[(c)]
The knots $5_2$, $6_3$, the Hopf link ($2_{1}^{2}$), the Borromean rings 
($6_{2}^{3}$) and the unknot ($T_1$) are $t_3$ equivalent.

\item[(d)]
The figure eight knot ($4_1$) and the knot $9_{42}$
(in the Rolfsen notation \cite{Rol}) are $t_3$ equivalent.

\item[(e)]
No links from different classes ((a), (b),
(c), (d)) are $t_3$ equivalent however links of (c) and  (d) 
 are  $t_3$, $\bar t_3$ equivalent (i.e. there is a sequence of  
$t_{3}^{\mp 1}$ or  $\bar t_{3}^{\mp 1}$ moves which lead from one link
 to another) and there is no more $t_3$, $\bar t_3$ equivalences among the above links.
\end{enumerate}

The $t_3$ and  $t_3$, $\bar t_3$ equivalences are illustrated in Fig. 3.10, 3.11, 3.12 and 3.13.

\vspace*{0.8in} \centerline{\psfig{figure=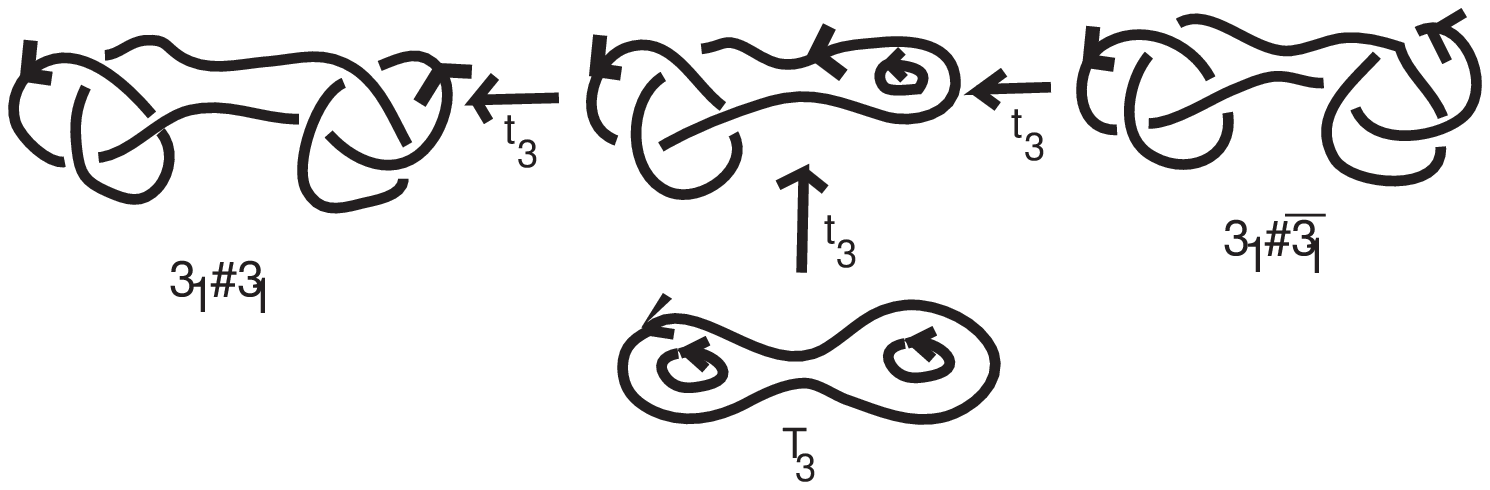,height=4.5cm}}
\begin{center}
Fig. 3.10
\end{center}

\vspace*{0.8in} \centerline{\psfig{figure=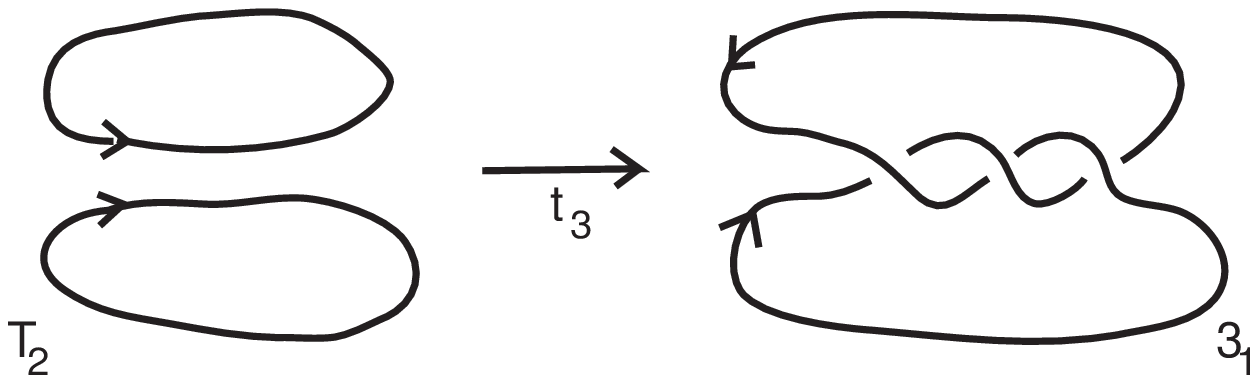,height=4cm}}
\begin{center}
Fig. 3.11
\end{center}

\vspace*{0.8in} \centerline{\psfig{figure=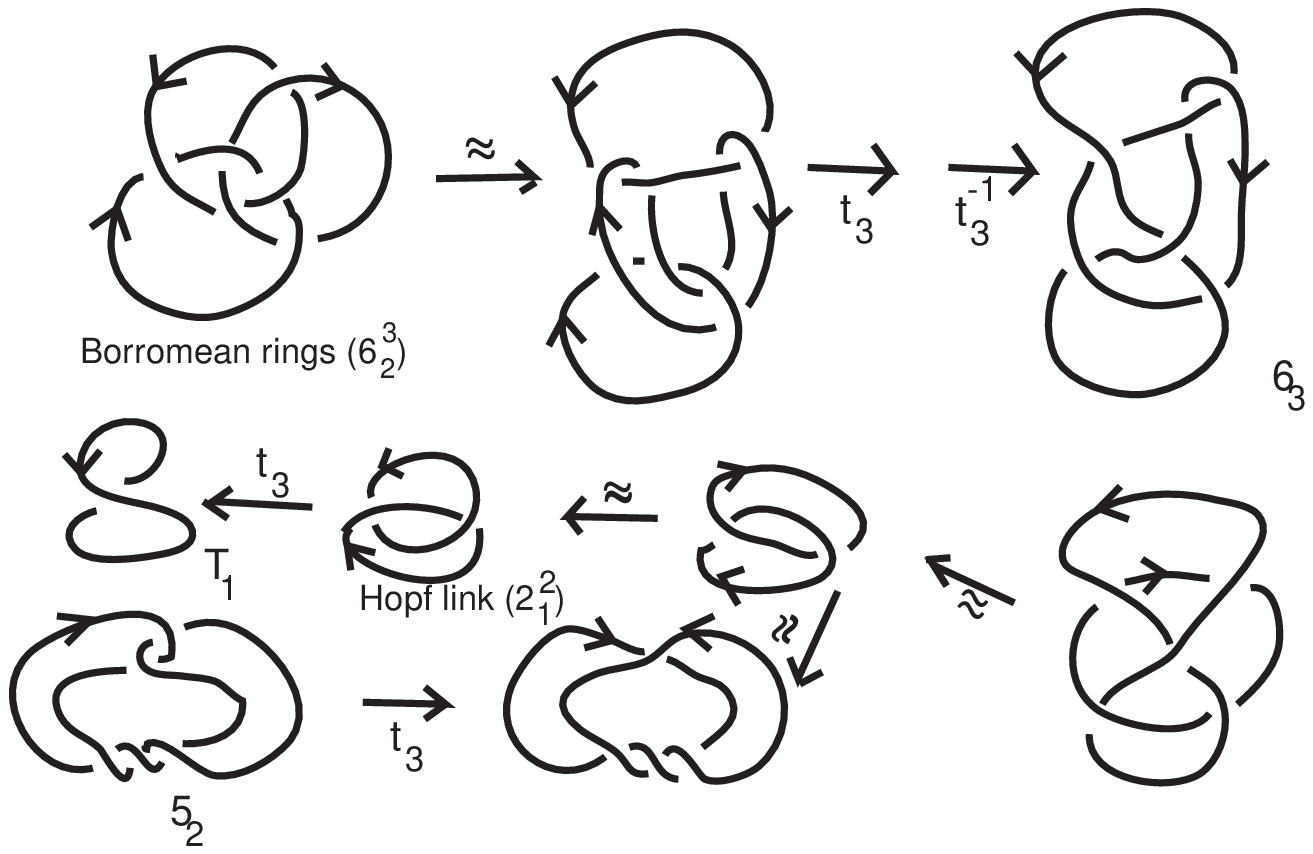,height=9cm}}
\begin{center}
Fig. 3.12
\end{center}

\vspace*{0.8in} \centerline{\psfig{figure=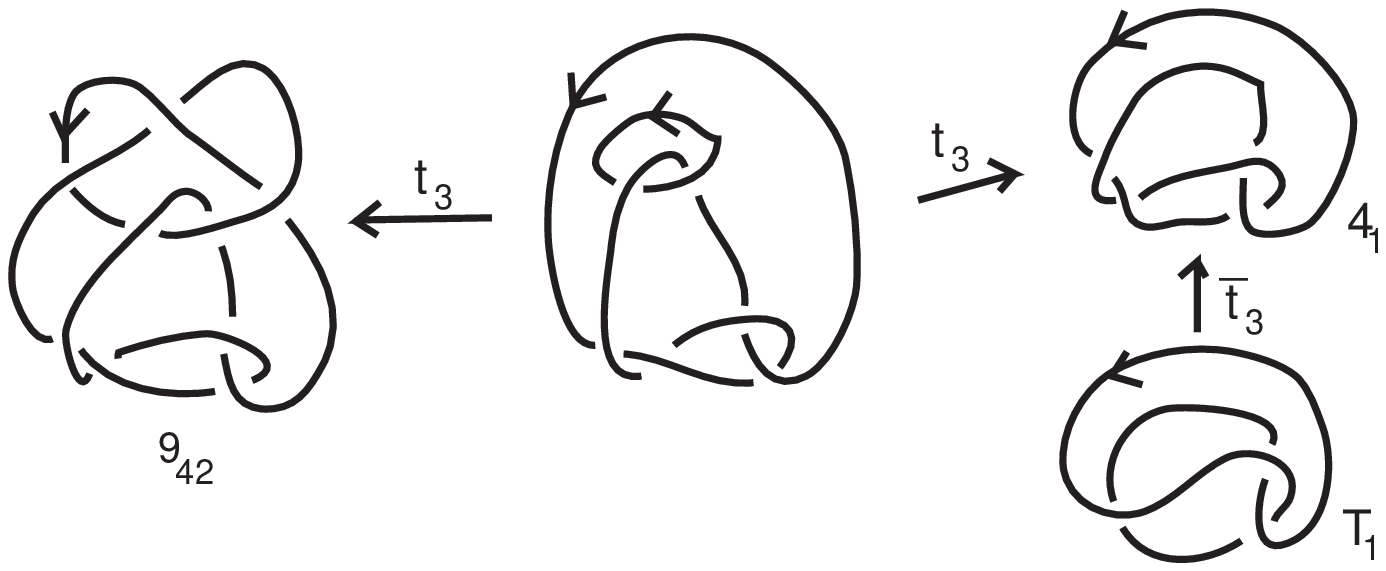,height=4.5cm}}
\begin{center}
Fig. 3.13
\end{center}

The first part of (e) follows Corollary 1.2 (a $t_3$ move changes $P_L(a,1)$ 
by the factor $-a^{-3}$) and the following computation:

$$P_{T_1}(a,1)=1,\ P_{T_2}(a,1)=a+a^{-1},\
P_{T_3}(a,1)=a^{-2}+2+a^2,\ P_{4_1}(a,1)=-a^{-2}-a^2.$$

The last statement of (e) follows from the fact that different trivial links are
 not $t_3$, ${\bar{t}}_3$ equivalent (see  Lemma 3.10(c) below).
\end{example}

\vspace{20pt}

\begin{lemma}
Consider the Jones polynomial $V_L(t)$ for $t=e^{\pi i/3}$ ($t^{1/2}=-e^{\pi i/6}$), then

\begin{enumerate}
\item[(a)]
$V_{t_3(L)}(e^{\pi i/3})=iV_L(t)$

\item[(b)]
\vspace{1mm}
   \renewcommand{\arraystretch}{2}
   $V_{{\bar{t}}_3(L)}(e^{\pi i/3})=\left\{
   \begin{array}{lr}
  {\rm(-1)}^{\lambda}iV_L(t) &
   {\rm if~ two~ components~ of}~L~{\rm\ are~ involved~ in}~{\bar{t}}_3~ {\rm{move}}\\
  {\rm(-1)}^{\lambda}V_L(t) &
   {\rm if~ one~ component~ of}~L~{\rm\ is~ involved~ in}~{\bar{t}}_3~{\rm move}\\
   \end{array}
   \right. $
   \par
   \vspace{2mm}

$\lambda$ depends on the linking numbers of components of $L$ and ${\bar{t}}_3(L)$
 and on an orientation of ${\bar{t}}_3(L)$ (see Theorem 1.13).

\item[(c)]
The trivial links $T_k$ and $T_j$ ($k\neq j$) are not $t_3$, 
${\bar{t}}_3$ equivalent and  $V_{T_k}(e^{\pi
i/3})=(\sqrt{3})^{k-1}$.
\end{enumerate}
\end{lemma}

\begin{proof}
(a) follows from Corollary 1.4, and (b) from Theorem 1.13. (c) follows from  (a) and (b).
\end{proof}

\vspace{20pt}

\begin{example}
Every closed $3$-braid link is $t_3$ equivalent 
to $T_1$, $T_2$, $T_3$ or the figure eight knot 
(in fact one can go from any closed $3$-braid link to one of these links using 
$t_3$ moves and regular isotopy). 
Because all presentations of $4_1$ and $T_3$ as closed $3$-braids have the same exponent sum (equal to $0$) 
therefore for any link $L$ which is $t_3$ equivalent to $4_1$ or 
$T_3$, each of its presentation as a $3$-braid has the same exponent sum  
(equal to  $3|4_1,L|_{t_3}^{\rm
lev}$ or $3|T_3,L|_{t_3}^{\rm lev}$). 
Consider, for example, the closed $3$-braid knot 
$\delta_{1}^{4}\delta_{2}^{-1}\delta_1\delta_{2}^{-4}$ (Fig. 
3.14). It is  $t_3$ equivalent to the figure eight knot so now we know that all presentations of this knot as a 
$3$-braid have the exponent sum equal to zero; on the other hand,
 the knot is  $t_4$ equivalent to the unknot so the method of Example 3.4  would not suffice to get 
the unique exponent sum.

\vspace*{0.8in} \centerline{\psfig{figure=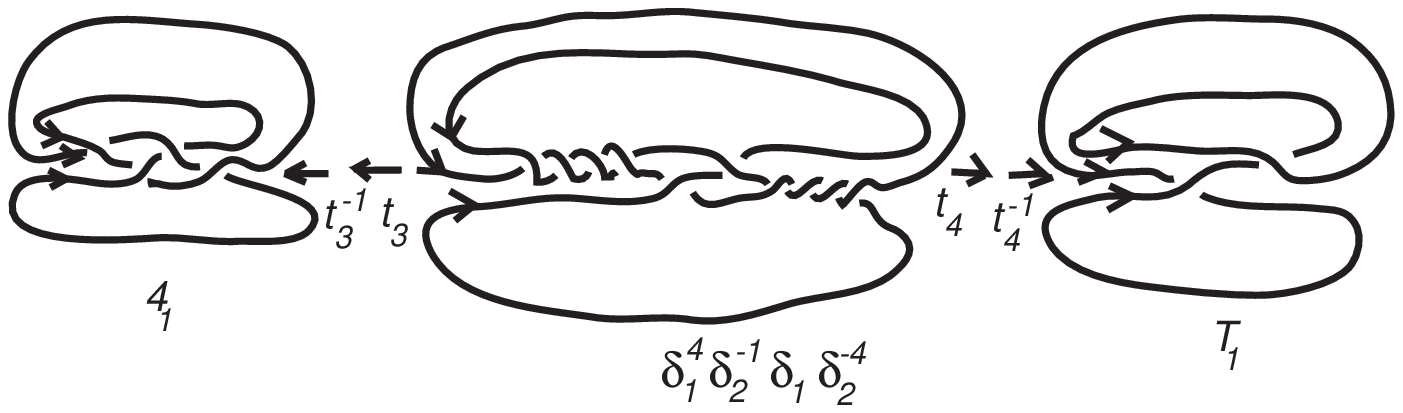,height=5cm}}
\begin{center}
Fig. 3.14
\end{center}

\end{example}

\vspace{20pt}

\begin{example}
Consider the following theorem of W.B.R.Lickorish and  
K.Millett \cite{L-M-2} (conjectured by J.Birman and partially proved by V.Jones). 

$$P_L(e^{\pi i/6},1)=V_L((e^{\pi i/3})=\mp
i^{c(L)-1}(i\sqrt{3})^{{\rm Dim\ }H_1(M_{L}^{(2)},{\bf Z}_3)},$$

where  $c(L)$ denotes the number of components of  $L$ and $t^{1/2}=e^{-\pi i/6}$ in  $V_L(t)$.

It follows from Lemma 3.10 that $t_3$ and $\bar t_3$ moves change  $V_L((e^{\pi i/3})$
 by factors $\mp 1$ or $\mp i$ and the second case happens if the move changes the number
of components. On the other hand $t_3$ and $\bar t_3$ moves 
preserve  $H_1(M_{L}^{(2)},{\bf Z}_3)$ (Theorem 3.1(b)) 
and for the trivial link  $T_n$, ${\rm Dim\
}H_1(M_{L}^{(2)},{\bf Z}_3)=n-1$.  Therefore the formula of 
Lickorish-Millett holds immediately from the above observations for a link which is 
$t_3$, $\bar t_3$ equivalent to a trivial link (the  sign  in formula can be found using Lemma 3.10;
 it was identified generally by  A.Lipson \cite{Lip}). This should be confronted with the following conjecture.

\vspace{20pt}

\begin{conjecture}(Montesinos-Nakanishi).
Every link is   $t_3$, $\bar t_3$ equivalent to a trivial 
link.\footnote{Added for e-print:\ 
The conjecture has been disproved in \cite{D-P-1}. 
The smallest known counter-example has 20 crossings.}
\end{conjecture}

\vspace{20pt}

It is an easy (but tedious) task to check the conjecture for closed $n$-braids 
($n\leq 5$) and  $n$-bridge links ($n\leq 3$) because for the braid group 
 $B_n$ ($n\leq 5$) the group  $B_n/(\delta_{1}^{3})$ is finite (\cite{Cox}),
 however the author did it only for closed $3$-braids and 
$2$-bridge links.\footnote{Added for e-print:\ The conjecture holds for 
 4-bridge links \cite{P-Ts,Tsu}. Furthermore every closed 5-braid is 
 $t_3$, $\bar t_3$ equivalent to a trivial link or to the closure of the 
5-string braid $(\delta_1\delta_2\delta_3\delta_4)^{10}$ \cite{Chen}. 
The last link is a counter-example to Montesinos-Nakanishi conjecture 
\cite{D-P-1}.}

\end{example}

\vspace{20pt}

\begin{example}
Consider the following theorem of Lickorish and 
Millett \cite{L-M-2} and H. Murakami \cite{Mur-2}:

$$P_L(1,1)=(2)^{(1/2){\rm Dim\ }H_1(M_L^{(3)},{\bf Z}_2)}.$$

It follows from Corollaries 1.2 and 1.8 that  
$t_3$ and $\bar t_4$ moves preserve 
 $P_L(1,1)$. Furthermore, $P_{T_n}(1,1)=2^{n-1}$. On the other hand ${\rm Dim\ }H_1(M_{T_n}^{(3)},{\bf Z}_2)=2(n-1)$ and 
$\bar t_4$ moves preserve  $H_1(M_{L}^{(3)},{\bf Z}_2)$. 
It can be shown, using the Fox approach that  $t_3$-moves preserve
 $H_1(M_{L}^{(3)},{\bf Z}_2)$ (\cite{P-3}). Therefore the formula of 
Lickorish-Millett-Murakami follows immediately from the above observations for a link which is 
$t_3$, $\bar
t_4$ equivalent to a trivial link (i.e. a link which can be got from a trivial one using 
$t_3$ and  $\bar t_4$ moves and isotopy). This leads to the following conjecture.

\vspace{20pt}

\begin{conjecture}
Every link is   $t_3$, $\bar t_4$ equivalent to a trivial link.
\end{conjecture}

\vspace{20pt} The author has checked the conjecture for closed $3$-braid links 
(see the remark after Conjecture 3.13).\footnote{Added for e-print:\ It has 
been checked for closed $4$-braid links \cite{Chen}.}

For  $t_5$ and  $\bar
t_4$ moves the analogy of Conjecture 3.15 does not hold.
 For example, trivial links, the trefoil knot, $8_5$ knot and  $8_{18}$ knot (\cite{Rol}, see 
Fig. 3.15) are not pairwise  $t_5$, $\bar t_4$ equivalent.
 The reason is that by Corollary 1.8 a $\bar t_4$ move does not change 
$P_L(1,z)$ and by Corollary 1.2 $t_5$ move changes 
$P_L(1,\frac{1+\sqrt{5}}{2})$ by the factor  $-1$ (notice that 
$2\cos (\pi i/5)=\frac{1+\sqrt{5}}{2}$); on the other hand all mentioned 
above links have pairwise different absolute values of $P_L(1,\frac{1+\sqrt{5}}{2})$:

$$P_{T_n}(1,\frac{1+\sqrt{5}}{2})=(\sqrt{5}-1)^{n-1}$$

$$P_{3_1}(1,\frac{1+\sqrt{5}}{2})=\frac{-3+\sqrt{5}}{2},\
P_{8_5}(1,\frac{1+\sqrt{5}}{2})=-4+\sqrt{5},$$

$$P_{8_{18}}(1,\frac{1+\sqrt{5}}{2})=\frac{1-2\sqrt{5}}{2}.$$

\vspace*{0.7in} \centerline{\psfig{figure=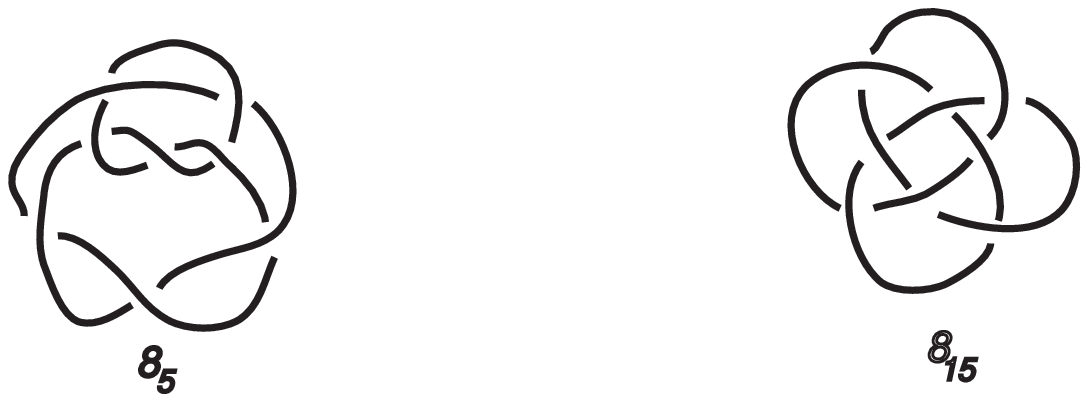,height=3.9cm}}
\begin{center}
Fig. 3.15
\end{center}

\end{example}

There is no chance for anything analogous to Conjectures 3.6 or 3.13 for  $t_k$, 
$\bar t_k$ moves, $k\geq 5$ (i.e. that all links are $t_k$, $\bar t_k$ equivalent to the trivial links).
 In particular  V.Jones (\cite{Jo-3}; Corollary  14.7) proved that the set 
$\{|V_L(e^{\pi i/5})|:\ L {\rm is\ a\ link}\}$ is dense in $\langle 0,\infty )$.  
On the other hand  $t_5$ and  $\bar
t_5$ moves do not change the absolute value of $V_L(e^{\pi i/5})$ (see Corollary  1.4 and  Theorem 1.13), 
and for trivial links, the values  $|V_{T_n}(e^{\pi
i/5})|=(2\cos \pi/10)^{n-1}$ are discrete in $\langle 0,\infty
)$.

There are natural relations between $t_k$ moves and signatures of links; 
we will list here some examples of such relations. For convenience, we start from 
the definition of the Tristram-Levine signature  (see  \cite{Gor,P-T-2} or 
\cite{P-1}). Let  $A_L$ be a Seifert matrix of a link $L$. 
For each complex number  $\xi$ ($\xi\neq 1$) 
consider Hermitian matrix  $A_L(\xi)=(1-\bar\xi
)A_L+(1-\xi )A_{L}^{T}$. The signature of this matrix, 
$\sigma_L(\xi )$  is called the 
Tristram-Levine signature of the link $L$.
 The classical signature  $\sigma$ 
satisfies  $\sigma_L=\sigma_L(0)$.

\vspace{20pt}

%%%%%%%%%%%%%%%%%%%%%%%%'±'±'©'ç

\begin{theorem}
\begin{enumerate}
\item[(a)]
For any $t_k$ move on an oriented link $L$

$$k-2\leq\sigma_L -\sigma_{t_k(L)}\leq k,$$

\item[(b)]
$0\leq\sigma_{\bar t_{2k}(L)}-\sigma_L(\xi)\leq 2$ if ${\rm
Re} (1-\xi )\geq 0$,

\item[(c)]
$\sigma_L(\xi_0)-\sigma_{t_4}(\xi_0)=2$ if ~
$P(i,\sqrt{2})\neq 0$ ~and ~ $\xi_0=1-e^{\pi
i/4}=\frac{2-\sqrt{2}}{2}-i\frac{\sqrt{2}}{2}$,

\item[(d)]
If two links $L_1$~ and~ $L_2$ are  $t_4$ equivalent~ then~
$|L_1,L_2 |_{t_4}^{\rm
lev}=(1/2)(\sigma_{L_1}(\xi_0)-\sigma_{L_2}(\xi_0))$,
~ provided  $P_{L_1}(i,\sqrt{2})\neq 0$.

\end{enumerate}
\end{theorem}

\begin{proof}
(a) We use the formula of  C.McA.Gordon, R.A.Litherland \cite{G-L} and ~
A.Marin, which links the signature of  Goeritz ~matrix~
of a link with a classical signature. ~We use~ the same notation as in the 
proof of Theorem 3.1. 
Divide the crossings of a given oriented link $L$ into two types as shown in Fig.~3.16.

\vspace*{0.8in} \centerline{\psfig{figure=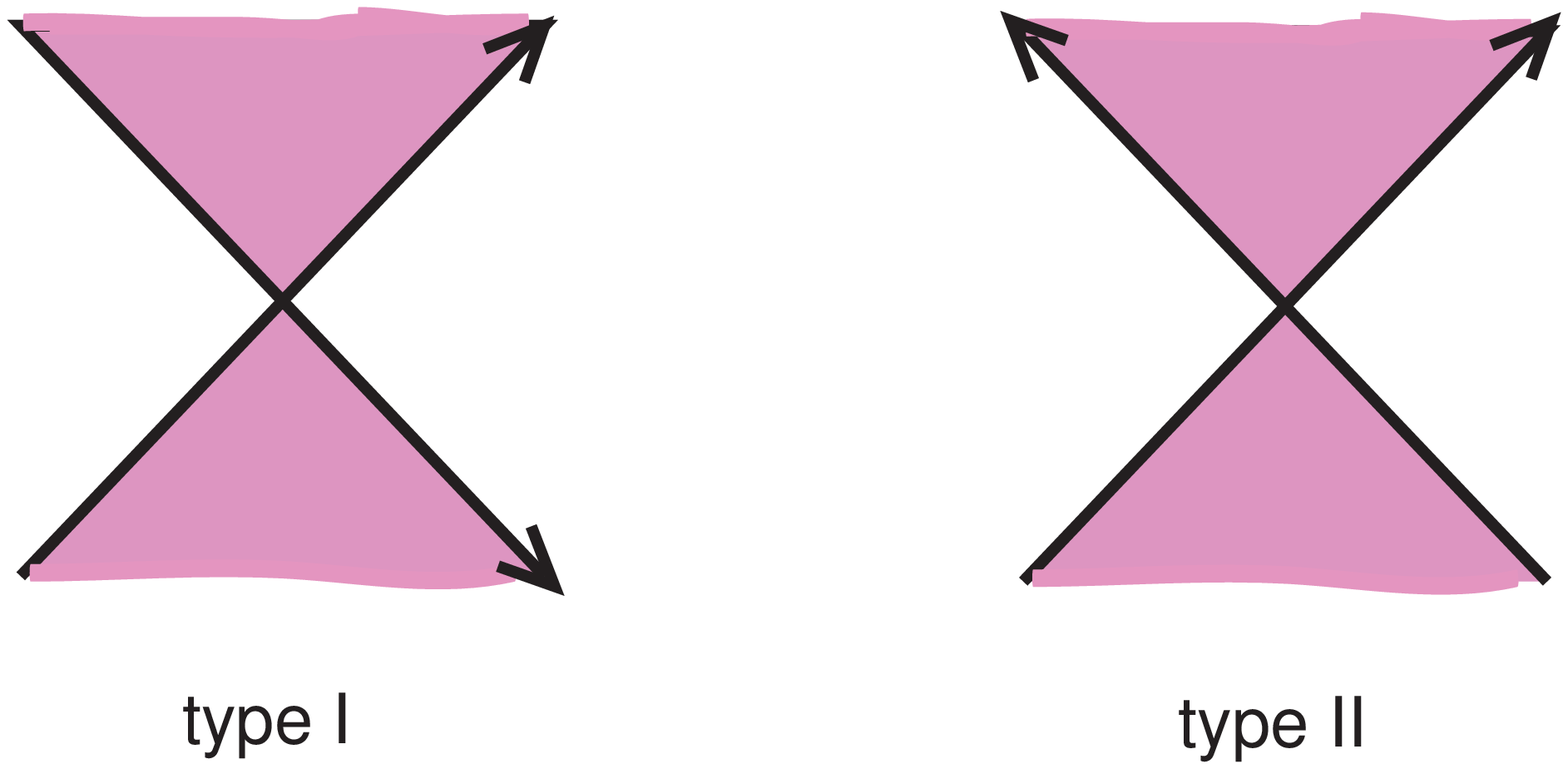,height=3cm}}
\begin{center}
Fig. 3.16
\end{center}

Define  $\mu=\sum\eta (p)$, summed over all crossing points of type II then  $\sigma_L=\sigma
(G_L)-\mu (L)$. Furthermore we have  $\mu (t_k(L))-\mu (L)=k$ (see Fig. 3.2), 
and from the form of the matrices  $G_{L'}=G_{t_k(L)}$ and  $G_L$ 
(see proof ~of Theorem 3.1) follows that  
$-2\leq \sigma(G_L)-\sigma(G_{t_k(L)})\leq 0$ and therefore 
$-2\leq \sigma_L+\mu
(L)-\sigma_{t_k(L)}-\mu (t_k(L))\leq 0$ and Theorem 3.16(a) follows.

To prove (b), we have to choose a proper  Seifert surface from 
which we will find the adequate 
Seifert matrix so one could easily compare the Levine-Tristram 
signature for $L$  and  $\bar t_{2k}(L)$. 
We can assume that Seifert surfaces for $L$ and $\bar
t_{2k}(L)$ looks locally as on Fig. 3.17 (or 
$\bar t_{2k}(L)$ s isotopic to  $L$).

\vspace*{0.8in} \centerline{\psfig{figure=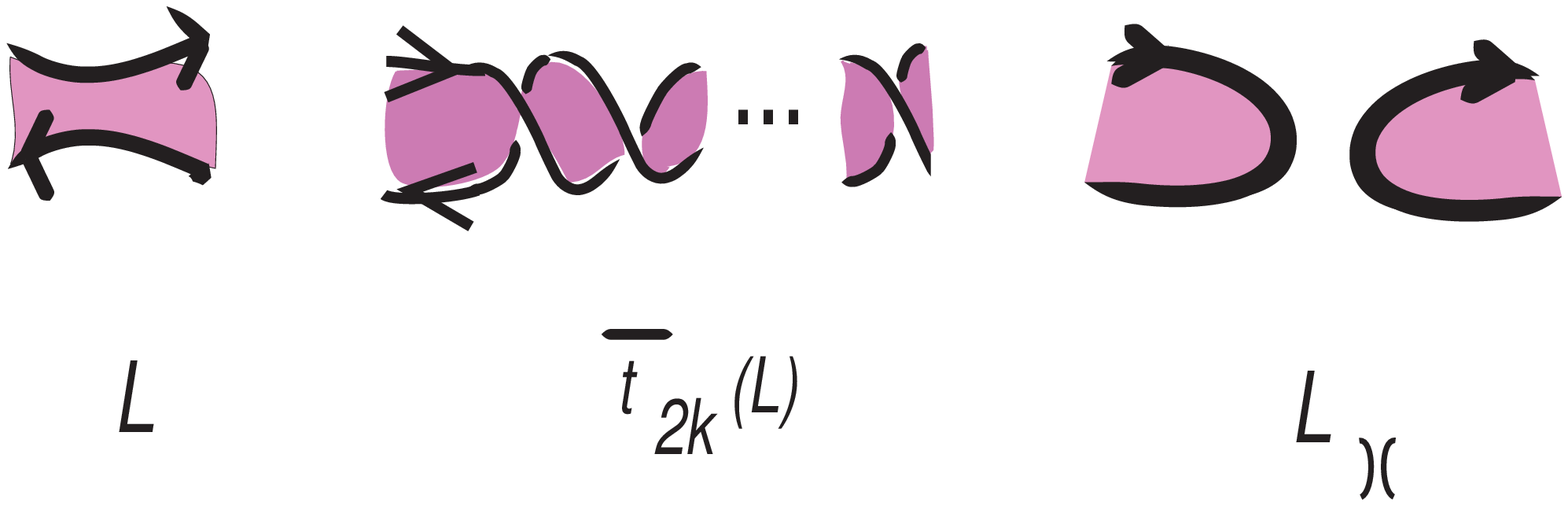,height=4.5cm}}
\begin{center}
Fig. 3.17
\end{center}

Then the Seifert matrices (in appropriate basis) are of the form :

  $A_{\bar t_{2k}(L)}=\left[
  \begin{array}{cc}
  A_{L_{\psfig{figure=inftyunoriented.eps,height=0.4cm}}} & \alpha \\
   \beta & q+k
   \end{array}
   \right],$

    $A_{L}=\left[
  \begin{array}{cc}
  A_{L_{\psfig{figure=inftyunoriented.eps,height=0.4cm}
   }} & \alpha \\
   \beta & q
   \end{array}
   \right],$

   where  $A_L$ is the Seifert matrix of $L$, $\alpha$ is a column, ~
   $\beta$ is a row and $q$ is a number  (compare \cite{Ka-1,P-T-2} or \cite{P-1}).
   Therefore 

  $A_{\bar t_{2k}(L}(\xi)=\left[
  \begin{array}{cc}
  A_{L_{\psfig{figure=inftyunoriented.eps,height=0.4cm}
   }}(\xi) & a \\
   a^{-T} & m+k(2-\xi-\bar\xi)
   \end{array}
   \right],$

    $A_{L}(\xi)=\left[
  \begin{array}{cc}
  A_{L_{\psfig{figure=inftyunoriented.eps,height=0.4cm}
   }}(\xi) & a\\
   a^{-T} & m
   \end{array}
   \right],$

   where  $a=(1-\bar\xi)\alpha+(1-\xi)\beta^T$ and 
   $m=((1-\bar\xi)+(1-\xi))q$. Because  $2-\xi-\bar\xi\geq 0$, 
   so 
   $0\leq\sigma (A_{\bar t_{2k}(L)}(\xi))-\sigma(A_{L}(\xi))\leq
   2$~
   and the proof of  (b) is finished.

   To prove  (c) we need further characterization of the Tristram-Levine signature, 
   
   given in \cite{P-T-2} (see also~
   \cite{P-1}); Assume  $|1-\xi |=1$, we have :~

   \begin{enumerate}
\item[(i)]
${\rm Det\ }iA_L(\xi )=P_L(i,2-\xi-\bar\xi)=\Delta_L(t')$ (for ~
$\sqrt{t'}=-i(1-\xi )$),

\item[(ii)]
$i^{\sigma_L(\xi)}=\frac{\Delta_L(t')}{|\Delta_L(t')|}$ if~
$\Delta_L(t')\neq 0$,

\item[(iii)]
$0\leq \sigma_L(\xi)-\sigma_{t_4(L)}(\xi)\leq 4$ ~if~ ${\rm
Re}(1-\xi)\geq 0$.
\end{enumerate}

((iii) can be got using  (b) two times with $k=1$; $\bar
t_2$ moves are equivalent to  $t_2$ moves).

Now consider the case when  
$2-\xi_0-\bar\xi_0=\sqrt{2}$ ($1-\xi_0=e^{\pi i/4}$). Then by Corollary 1.3,~
$\Delta_{t_4(L)}(t')=-\Delta_L(t')$. Therefore by  (ii) and (iii)~
 $\delta_L(\xi_0)-\delta_{t_4(L)}(\xi_0)=2$. (d) follows immediately from (c).~
\end{proof}

\vspace{20pt}

One can expect interesting relations between $t_k$ moves  and non-cyclic coverings of links.~
We limit ourself to two examples, first of which was suggested by R.Campbell. 

\vspace{20pt}

\begin{example}
\begin{enumerate}
\item[(a)]
A link diagram is $3$-coloured if every overpass is coloured,~
say, red, yellow or blue, at least two coloures are used and 
at any given crossing either all three colours appear or only one colour 
appears \cite{Fo-2}. 
Then if a link $L_1$ is  $t_3$, $\bar t_3$ equivalent to  $L_2$ then 
either both links are $3$-coloured or none of them are $3$-coloured. In particular a link which is  $t_3$, $\bar t_3$ equivalent to a trivial link of more than one component is 
$3$-coloured.~
The proof is illustrated in Fig. 3.18. The link  $6_{3}^{2}$~
\cite{Rol} is $3$-coloured in Fig. 3.19.

\vspace*{0.8in} \centerline{\psfig{figure=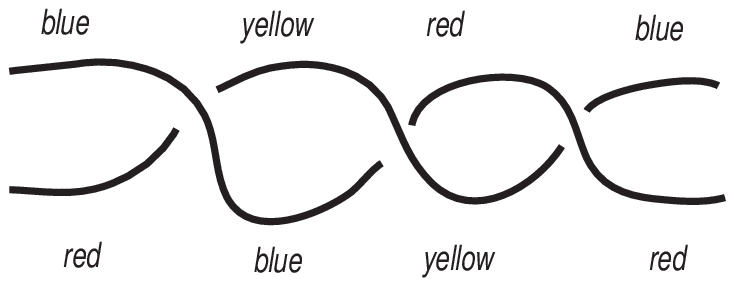,height=4.5cm}}
\begin{center}
Fig. 3.18
\end{center}

\item[(b)]
$3$-colouring corresponds to an epimorphism  $\pi_1(S^3-L)\rightarrow
S_3$; more generally we have:
If a knot $K_1$ is  $t_{2p}$, $\bar t_{2p}$ equivalent to  $K_2$ ($p$ -
prime) then  either both knots or none of them have dihedral representations 
i.e. epimorphism~

$$\pi_1(S^3-K)\rightarrow D_{2p}=\{a,b:a^2=1, b^p=1, a b
a=b^{-1}\}.$$

It follows from the fact that  $t_{2p}$, $\bar t_{2p}$ moves preserve 
$H_1(M_{K}^{(2)}, {\bf Z}_p)$ (Theorem 3.1(b)) and from the result 
of Fox that the epimorphism exists iff 
 $H_1(M_{K}^{(2)}, {\bf Z}_p)$, 
is nontrivial \cite{Fo-2} (see also  \cite{B-Z}; 14.8).\footnote{Added 
for e-print:\ These ideas have been developed in \cite{P-5}.}

\vspace*{0.6in} \centerline{\psfig{figure=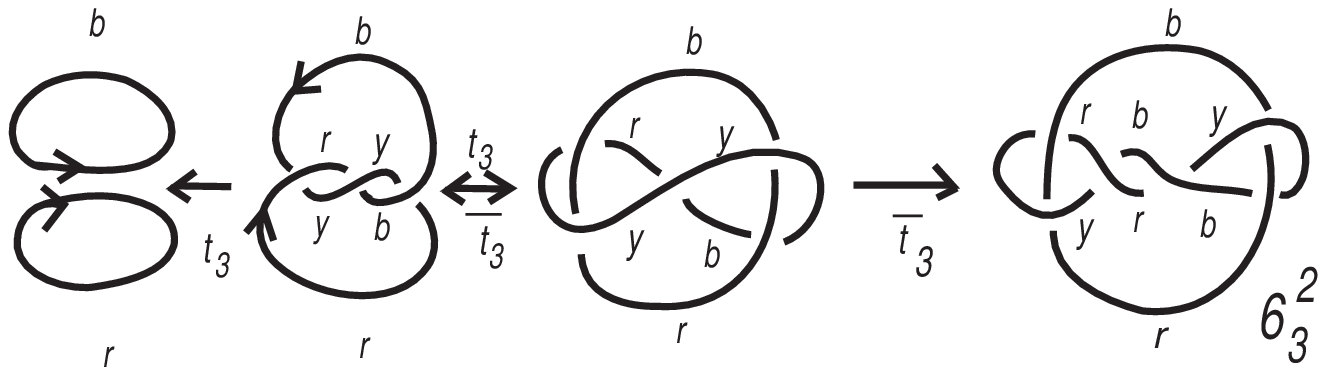,height=5.5cm}}
\begin{center}
Fig. 3.19
\end{center}

\end{enumerate}
\end{example}

\vspace{30pt}
%\newpage

\vspace{20pt}

\centerline{Department of Mathematics}

\centerline{Warsaw University}

\centerline{00901 Warszawa, PKiN IXp.}

\centerline{Poland}


\begin{thebibliography}{99}

\bibitem [A]{A}
J.W.Alexander, {\em Topological invariants of knots and links},
Trans. Amer. Math. Soc., 30 (1928), 275-306.

\bibitem [Bi]{Bi}
J.S.Birman, {\em Braids, links and mapping class groups}, Ann.
Math. Studies 82, Princeton Univ. Press, 1974.

\bibitem [B-W]{B-W}
J.S.Birman, B.Wajnryb, {\em Markov classes in certain finite
quotients of Artin's braid group}, Israel J. Math. (to appear).
Added for e-print:\ 56(2), 1986, 160-178.

\bibitem [B-Z]{B-Z}
G.Burde, H.Zieschang, {\em Knots}, De Gruyter Studies in Math. 5,
Berlin, New York 1985.

\bibitem[Chen]{Chen}
Added for e-print:\ 
Q. Chen, The $3$-move conjecture for $5$-braids,
Knots in Hellas' 98;
The Proceedings of the International Conference
on Knot Theory and its Ramifications; Volume 1.
In the Series on Knots and Everything, Vol. 24, September 2000,
pp. 36-47.

\bibitem [Con]{Con}
J.H.Conway, {\em An enumeration of knots and links, and some of
their algebraic properties}, Computational problems in abstract
algebra (J.Leech, ed.), Pergamon, Oxford and New York (1969),
329-358.

\bibitem [Cox]{Cox}
H.S.M.Coxeter, {\em Factor groups of the braid group}, Proc.
Fourth Canadian Math. Congress, Banff, 1957, 95-122.

\bibitem [D-P-1]{D-P-1} Added for e-print:\
M.~K.~D{\c a}bkowski, J.~H.~Przytycki, Burnside
obstructions to the Montesinos-Nakanishi 3-move conjecture, i
{\it Geometry and Topology (G\&T)}, 6, June, 2002, 335-360;\\
e-print:\ http://front.math.ucdavis.edu/math.GT/0205040

\bibitem [D-P-2]{D-P-2} Added for e-print:\ 
M.~K.~D{\c a}bkowski, J.~H.~Przytycki, Unexpected
connection between knot theory and Burnside groups,
{\it Proc. Nat. Acad. Science},
101(50), December, 2004, 17357-17360; \\
 e-print: \ http://front.math.ucdavis.edu/math.GT/0309140

\bibitem [Fo-1]{Fo-1}
R.H.Fox, {\em Congruence classes of knots}, Osaka Math. J. 10
(1958), 37-41.

\bibitem [Fo-2]{Fo-2}
R.H.Fox, {\em Metacyclic invariants of knots and links}, Canadian
J. Math. XXII (2) (1970), 193-201.

\bibitem[F-Y-H-L-M-O]{F-Y-H-L-M-O}
P.Freyd, D.Yetter, J.Hoste, W.B.R.Lickorish, K.Millet, A.Ocneau,
{\em A new polynomial invariant of knots and links}, Bull. Amer.
Math. Soc. 12(2) (1985), 239-249.

\bibitem[Goe]{Goe}
L.Goeritz, {\em Knoten und quadratische Formen}, Math. Z. 36
(1933), 647-654.

\bibitem [Gor]{Gor}
C.McA.Gordon, {\em Some aspects of classical knot theory}, In: Knot
theory, L.N.M. 685 (1978), 1-60.

\bibitem [G-L]{G-L}
C.McA.Gordon, R.A.Litherland, {\em On the signature of a link},
Inv. Math. 47 (1978), 53-69.

\bibitem [Jo-1]{Jo-1}
V.F.R.Jones, {\em A polynomial invariant for knots via von Neumann
algebras},  Bull. Amer. Math. Soc. 12(1) (1985), 103-111.


\bibitem [Jo-2]{Jo-2}
V.F.R.Jones, {\em A new knot polynomial and von Neumann algebras},
Notices AMS 33(2) (1986), 219-225.

\bibitem [Jo-3]{Jo-3}
V.F.R.Jones, {\em Hecke algebra representations of braid groups
and link polynomials}, Ann. of Math. (to appear).
Added for e-print:\ 126(2), 1987, 335-388. 

\bibitem [Kan]{Kan}
T.Kanenobu, {\em Examples on polynomial invariants of knots and
links}, Math. Ann. 275 (1986), 555-572.

\bibitem [Ka-1]{Ka-1}
L.H.Kauffman, {\em The Conway polynomial}, Topology 20 (1980),
101-108.

\bibitem [Ka-2]{Ka-2}
L.H.Kauffman, {\em Knots}, Lecture notes, Zaragoza, Spring 1984.
\ \ Added for e-print:\
{\it On knots}, Annals of Math. Studies, 115,
Princeton University Press, 1987.

\bibitem [Ka-3]{Ka-3}
L.H.Kauffman, {\em An invariant of regular isotopy}, preprint
1985.\ \ Added for e-print:\ 
{\it Trans. Amer. Math. Soc.}, 318(2), 1990, 417--471.

\bibitem [Ka-4]{Ka-4}
L.H.Kauffman, {\em State models for knot polynomials}, preprint
1985. \ \ Added for e-print:\ 
{\it Topology} 26, 1987, 395-407.

\bibitem [Ki]{Ki}
M.E.Kidwell, {\em Relations between the Alexander polynomial and
summit power of a closed braid}, Math. Sem. Notes Kobe 10 (1982) 387-409.

\bibitem [Kin-1]{Kin-1}
S.Kinoshita, {\em On Wendt's theorem of knots, I}, Osaka Math. J.
9(1) (1957) 61-66.

\bibitem [Kin-2]{Kin-2}
S.Kinoshita, {\em On Wendt's theorem of knots, II}, Osaka Math. J.
10 (1958) 259-261.

\bibitem [Kin-3]{Kin-3}
S.Kinoshita, {\em On the distribution of Alexander polynomials of
alternating knots and links}, Proc. Amer. Math. Soc. 79(4) (1980),
644-648.

\bibitem[Li]{Li}
W.B.R.Lickorish, {\em A relationship between link polynomials},
Math. Proc. Cambridge Phil. Soc. (1986), 100, 109-112.

\bibitem [L-M-1]{L-M-1}
W.B.R.Lickorish, K.C.Millet, {\em A polynomial invariant of
oriented links}, Topology 26(1) (1987), 107-141.

\bibitem [L-M-2]{L-M-2}
W.B.R.Lickorish, K.C.Millet, {\em Some evaluations of link
polynomials}, Com. Math. Helvetici 61(3) (1986), 349-359.

\bibitem [L-M-3]{L-M-3}
W.B.R.Lickorish, K.C.Millet, {\em The reversing result for the
Jones polynomial}, Pacific J. Math. (to appear). 
 \ \ Added for e-print:\
124 (1986), 173-176.

\bibitem [Lip]{Lip}
A.S.Lipson {\em An evaluation of a link polynomial}, Math. Proc.
Camb. Phil. Soc. 100 (1986), 361-364.

\bibitem [Mon]{Mon}
J.M.Montesinos {\em Lectures on 3-fold simple coverings and
3-manifolds}, Amer. Math. Soc. Contemp. Math. 44 (1985), 157-177.

\bibitem [Mo]{Mo}
H.R.Morton Personal conversation,  (July 1986).

\bibitem [Mu]{Mu}
K.Murasugi, {\em On closed 3-braids}, Memoirs AMS 151 (1974) Amer.
Math. Soc. Providence, R.I.

\bibitem [Mur-1]{Mur-1}
H.Murakami, {\rm A recursive calculation of the Arf invariant of
a link}, J. Math. Soc. Japan 38(2) (1986), 335-338.

\bibitem [Mur-2]{Mur-2}
H.Murakami, {\rm Unknotting number and polynomial invariants of a
link}, preprint 1985.   


\bibitem [N]{N}
Y.Nakanishi, {\em Fox's congruence classes and Conway's potential
functions of knots and links}, preprint 1986.

\bibitem [N-S]{N-S}
Y.Nakanishi, S.Suzuki, {\em On Fox's congruence classes of knots},
Osaka J. M. 24 (1987) (to appear).
 \ \ Added for e-print:\ pages 217--225.

\bibitem [P-1]{P-1}
J.H.Przytycki, {\em Survey on recent invariants in classical knot
theory}, preprint, Warsaw University, 1986.
 \ \ Added for e-print:\ Part of: {\it Teoria 
w\c ez\l\'ow: podej\'scie kombinatoryczne},
(Knots: combinatorial approach to the knot theory),
Script, Warsaw, August 1995, 240+XLVIIIpp. (Extended version 
for Cambridge University Press, to appear).

\bibitem [P-2]{P-2}
J.H.Przytycki, {\em $t_k$ equivalence of links and Conway formulas
for the Jones-Conway and Kauffman polynomials}, preprint 1986.
\ \ Added for e-print:\ 
{\it Bull. Polish Acad. Sci. Math.}, 36(11-12), 1988, 675-680.

\bibitem [P-3]{P-3}
J.H.Przytycki, {\em Plans' theorem for links: An application of
$t_k$ moves}, Bull. Canad. Math. Soc., (to appear).
\ \ Added for e-print:\ 31(3), 1988, 325-327.

\bibitem [P-4]{P-4}
J.H.Przytycki, {\em On Murasugi and Traczyk criteria for periodic
links}, preprint, March 1987.
\ \ Added for e-print:\ {\it Math. Ann.}, 283, 1989, 465 - 478.

\bibitem [P-5]{P-5} 
Added for e-print:\ 
J.~H.~Przytycki, 3-coloring and other elementary
invariants of knots, Banach Center Publications, Vol. 42,
{\it Knot Theory}, 1998, 275-295. 

\bibitem[P-T-1]{P-T-1}
J.H.Przytycki, P.Traczyk, {\em Invariants of links of Conway
type}, Kobe J. Math. (to appear).
Added for e-print:\  4, 1987, 115-139.


\bibitem[P-T-2]{P-T-2}
J.H.Przytycki, P.Traczyk, {\em Conway algebras and skein
equivalence of links}, preprint, 1985.
Added for e-print:\ Part of the paper published in\ \ 
{\it Proc. Amer. Math. Soc.}, 100(4), 1987, 744-748.

\bibitem [P-Ts] {P-Ts}
Added for e-print:\ 
J.~H.~Przytycki, T.Tsukamoto, The fourth skein module and
the Montesinos-Nakanishi conjecture for 3-algebraic links,
{\it  J. Knot Theory Ramifications}, 10(7), 2001, 959--982.\\
http://front.math.ucdavis.edu/math.GT/0010282

\bibitem [Re]{Re}
K.Reidemeister, {\em Knotentheorie}, Ergebn, Math. Grenzgeb. Bd.
1; Berlin: Springer-Verlag, 1932.

\bibitem [Rob]{Rob}
R.A.Robertello, {\em An invariant of knot cobordism}, Comm. Pure
Appl. Math. 18 (1965), 543-555.

\bibitem [Rol]{Rol}
D.Rolfsen, {\em Knots and links, Publish or Perish}, Inc. Berkeley
1976, Math. Lect. Series 7.

\bibitem [Tsu]{Tsu} 
Added for e-print:\ 
T.~Tsukamoto, The fourth skein module of 3-dimensional manifolds,
PhD dissertation, George Washington University, May, 2000.

\bibitem[Wa]{Wa}
B.Wajnryb, {\em Markov classes in certain finite symplectic
representations of braid groups}, these proceedings.
Added for e-print:\
Braids (Santa Cruz, CA, 1986), 687--695, Contemp. Math., 
78, Amer. Math. Soc., Providence, RI, 1988.

\bibitem [We]{We}
H.Wendt, {\em Die gordische Aufl\"{o}sung von Knoten}, Math. Z. 42
(1937), 680-696.
\end{thebibliography}
\end{document}